\documentclass{fundam}

\publyear{25}
\papernumber{4}
\volume{194}
\issue{3}
\theDOI{10.46298/fi.13125}

\versionForARXIV

\pdfoutput=1
\title{Dialectica Petri Nets}
\author{Elena Di Lavore\thanks{Elena Di Lavore was supported by the European Social Fund Estonian IT Academy research measure (project 2014-2020.4.05.19-0001).}\\
  Department of Software Science, Tallinn University of Technology, Tallinn, Estonia\\
  Department of Computer Science, University of Pisa, Pisa, Italy
  \and Wilmer Leal\thanks{Wilmer Leal was supported by the German Academic Exchange Service (DAAD): Forschungsstipendien-Promotionen in Deutschland, 2017/2018 (Bewerbung 57299294).}\\
  Department of Computer Science, University of Florida, Gainesville, Florida, USA\\
  Max Planck Institute for Mathematics in the Sciences, Leipzig, Germany \\
  Bioinformatics Group, Department of Computer Science, Universit\"at Leipzig, Leipzig, Germany
  \and Valeria de Paiva\\
  Topos Institute, Berkeley, California, USA}
\runninghead{E. Di Lavore, W. Leal, V. de Paiva}{Dialectica Petri Nets}
\usepackage{amssymb,amsmath,amsfonts,mathrsfs,graphicx,pdfpages,todonotes,bbold,subcaption}
\usepackage{cmll}

\usepackage{figures}

\newcommand{\Reals}{\mathbb{R}}
\newcommand{\Naturals}{\mathbb{N}}
\newcommand{\Integers}{\mathbb{Z}}

\newcommand{\id}[1]{\mathbb{1}_{#1}}
\newcommand{\swap}[1]{\sigma_{#1}}
\newcommand{\op}{^{op}} % opposite category
\newcommand{\Hom}[1]{\mathsf{Hom}_{#1}}
\newcommand{\diagonal}[1]{\Delta_{#1}}
\newcommand{\leftproarrow}{\mathrel{\ooalign{$\longleftarrow$\cr % arrow to the left with a bar
    \hidewidth\raisebox{.1ex}{$\shortmid$}\hidewidth\cr}}}
\newcommand{\dialarrow}[1]{\overset{#1}{\leftproarrow}}
\newcommand{\compose}{\circ}

\newcommand{\dialLineale}[1]{\mathsf{M}_{#1}\Set} % dialectica category on a lineale
\newcommand{\dialLinealeCat}[2]{\mathsf{M}_{#1}{\namedcat{#2}}} % dialectica category on a lineale
\newcommand{\Lin}{\mathsf{Lin}} % category of lineales
\newcommand{\dialNet}[1]{\mathsf{Net}_{#1}} % category of petri nets on a lineale
\newcommand{\product}{\mathbin{\&}}
\newcommand{\coproduct}{\oplus}
\newcommand{\productmap}[2]{\langle #1, #2\rangle}
\newcommand{\coproductmap}[2]{[#1,#2]}
\newcommand{\eval}{\mathsf{eval}}
\newcommand{\discard}{\epsilon} % discarding map in Set
\newcommand{\tensor}{\otimes}
\newcommand{\unit}{I}
\newcommand{\inthom}[2]{[#1,#2]}
\newcommand{\linorder}{\sqsupseteq}
\newcommand{\lintensor}{\ast}
\newcommand{\lininthom}{\multimap}
\newcommand{\tensortypesecond}[4]{#2^{#3} \times #4^{#1}} % type of the second component of the tensor product of objects
\newcommand{\inthomtypefirst}[4]{#3^{#1} \times #2^{#4}}

\newcommand{\pre}[1]{\prescript{\triangleright}{}{#1}}
\newcommand{\post}[1]{#1^{\triangleright}}

\newcommand{\namedcat}[1]{\mathsf{#1}}

\newcommand{\Cat}{\namedcat{Cat}}
\newcommand{\LinCat}{\namedcat{Lin}\Cat}
\newcommand{\DialCat}{\namedcat{Dial}\Cat}

\newcommand{\Set}{\namedcat{Set}}

\makeatletter
\newcommand*{\relrelbarsep}{.386ex}
\newcommand*{\relrelbar}{%
  \mathrel{%
    \mathpalette\@relrelbar\relrelbarsep
  }%
}
\newcommand*{\@relrelbar}[2]{%
  \raise#2\hbox to 0pt{$\m@th#1\relbar$\hss}%
  \lower#2\hbox{$\m@th#1\relbar$}%
}
\providecommand*{\rightrightarrowsfill@}{%
  \arrowfill@\relrelbar\relrelbar\rightrightarrows
}
\providecommand*{\leftleftarrowsfill@}{%
  \arrowfill@\leftleftarrows\relrelbar\relrelbar
}
\providecommand*{\xrightrightarrows}[2][]{%
  \ext@arrow 0359\rightrightarrowsfill@{#1}{#2}%
}
\providecommand*{\xleftleftarrows}[2][]{%
  \ext@arrow 3095\leftleftarrowsfill@{#1}{#2}%
}
\makeatother

\usepackage{subfiles}
\usepackage{stackengine}
\usepackage{mathtools}
\usepackage[utf8]{inputenc}
\usepackage{csquotes}
\usepackage{color}
\usepackage{tikz}
\usepackage{comment}

\usepackage{graphicx}
\usepackage{adjustbox}
\usepackage[all,2cell]{xy}\UseAllTwocells\SilentMatrices
\definecolor{darkgreen}{rgb}{0,0.45,0}
\usepackage{enumerate}

\usepackage{tikz}
\usepackage{tikz-cd}
\usetikzlibrary{backgrounds,circuits,circuits.ee.IEC,shapes,fit,matrix}
\tikzstyle{simple}=[-,line width=2.000]
\tikzstyle{arrow}=[-,postaction={decorate},decoration={markings,mark=at position .5 with {\arrow{>}}},line width=1.100]
\pgfdeclarelayer{edgelayer}
\pgfdeclarelayer{nodelayer}
\pgfsetlayers{edgelayer,nodelayer,main}

\tikzstyle{none}=[inner sep=0pt]

\definecolor{lblue}{rgb}{0,250,255}
\tikzstyle{species}=[circle,fill=yellow,draw=black,scale=1.15]
\tikzstyle{transition}=[rectangle,fill=lblue,draw=black,scale=1.15]
\tikzstyle{inarrow}=[->, >=stealth, shorten >=.03cm,line width=1.5]
\tikzstyle{empty}=[circle,fill=none, draw=none]
\tikzstyle{inputdot}=[circle,fill=black,draw=black, scale=.25]
\tikzstyle{inputarrow}=[->,draw=purple, shorten >=.05cm]
\tikzstyle{simple}=[-,draw=black,line width=1.000]

\usetikzlibrary{arrows,shapes,snakes,automata,backgrounds,petri}
\tikzstyle{place}=[circle,thick,draw=blue!75,fill=blue!20,minimum size=6mm]
\tikzstyle{red place}=[place,draw=red!75,fill=red!20]
\tikzstyle{transition}=[rectangle,thick,draw=black!75,
  			  fill=black!20,minimum size=4mm]

\tikzset{-|->/.style={decoration={markings,
      mark=at position 0.5 with {\arrow{|}},
      mark= at position 1 with{\arrow{>}}},
    postaction={decorate}}}

\definecolor{purple(x11)}{rgb}{0.8, 0, 0.8}

\definecolor{darkred}{rgb}{0.5, 0, 0}

\definecolor{penBlue}{HTML}{001E59}
\colorlet{localcolor}{penBlue}
\hypersetup{
  breaklinks=true,
  colorlinks=true,
  linkcolor=localcolor,
  citecolor=localcolor,
  urlcolor=black,
}

\begin{document}
\maketitle
\begin{abstract}
  The categorical modelling of Petri nets has received much attention recently. The Dialectica construction has also had its fair share of attention. We revisit the use of the Dialectica construction as a categorical model for Petri nets generalising the original application to  suggest that Petri nets with different kinds of transitions can be modelled in the same categorical framework. Transitions representing truth-values, probabilities, rates or multiplicities,  evaluated in different algebraic structures called lineales are useful and are modelled here in the same category. We investigate (categorical instances of)  this generalised model and its connections to more recent models of categorical nets.
\end{abstract}

\section{Introduction}

Petri nets exert endless fascination over category theorists.  Maybe category theorists see Petri nets as a gauntlet thrown at them, because the definition of a morphism of Petri nets is not obvious and different definitions lead to different categories.
Maybe the bipartite graphs that usually depict Petri nets look too  similar to automata ones, and these are  the initial  sources of good categorical examples in computing. In any case, many different categorical models of Petri nets do exist and some are fundamentally different from others. One fundamental difference is whether one concentrates on the token game and the behaviour of a given Petri net or on the graphs underlying different nets. Another difference is which possible operations combining different Petri nets one considers. Alternatively, the difference may lie in the type of relationships (labels) that the Petri net can model.

In this work, we present a categorical model of Petri nets with linear connectives. We explore the model originally introduced by Winskel~\cite{Winskel1987,winskel1988}, but use it with morphisms, as in the work of Brown~\cite{brown1990} and others,  that relate Petri nets to constructors in Linear Logic~\cite{girard1987}. These connectives enable the assembly of small networks into larger ones in a principled manner. The model is flexible enough to capture a broad range of relationships, provided that the set of labels encoding these relationships can be transformed into a lineale, a poset version of a symmetric monoidal closed category.

\paragraph{Petri nets.}

\begin{figure}[t]
  \begin{subfigure}[ht]{0.4\linewidth}
    \centering
    \includegraphics[width=.6\linewidth]{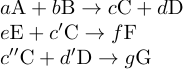}
  \end{subfigure}
  \hfill
  \begin{subfigure}[ht]{0.6\linewidth}
    \centering
    \begin{tikzpicture}[xscale=1.5]
      \node[place,tokens=0,label=center:A](A)at(0,4){};
      \node[place,tokens=0,label=center:B](B)at(0,2){};
      \node[place,tokens=0,label=center:C](C)at(2,4){};
      \node[place,tokens=0,label=center:D](D)at(2,2){};
      \node[place,tokens=0,label=center:E](E)at(2,6){};
      \node[place,tokens=0,label=center:F](F)at(4,5){};
      \node[place,tokens=0,label=center:G](G) at(4,3){};
      \node[transition,label=below:$r_1$](r1)at(1,3){}
      edge[pre]node[above]{$a$}(A)
      edge[post]node[above]{$c$}(C)
      edge[pre]node[above]{$b$}(B)
      edge[post]node[above]{$d$}(D);
      \node[transition,label=below:$r_2$](r2)at(3,5){}
      edge[post]node[above]{$f$}(F)
      edge[pre]node[above]{$c^{\prime}$}(C)
      edge[pre]node[above]{$e$}(E);
      \node[transition,label=below:$r_3$](r3)at(3,3){}
      edge[pre]node[above]{$c^{\prime\prime}$}(C)
      edge[pre]node[above]{$d^{\prime}$}(D)
      edge[post]node[above]{$g$}(G);
    \end{tikzpicture}
  \end{subfigure}%
  \captionsetup{width=.9\linewidth}
  \caption{Chemical reaction data represented as a list (left) and as a Petri net (right). $\text{A}, \text{B}, \ldots, \text{H}$ are substances and $a,b,\ldots, h$ are stoichiometric coefficients that indicate the proportion in which they combine.}\label{ChemDataRep}
\end{figure}

A \textit{Petri net} is simply a directed bipartite graph that has two types of elements, \textit{conditions} and \textit{events} (also called places and transitions). These are usually depicted as circles and rectangles, respectively (\Cref{ChemDataRep}, right).
This work is not concerned with the dynamic behaviour of Petri nets while it focuses on their structure.
However, keeping in mind the intended semantics of a net may help with intuitions.
When describing the dynamic behaviour of a Petri net, one starts with an \emph{initial marking}.
This consists of a number of \emph{tokens}, depicted as black circles inside places, for every place in the net as shown in \Cref{fig:initial-marking}.
\begin{figure}[h!]
  \centering
    \begin{tikzpicture}[xscale=1.5]
      \node[place,tokens=2](A)at(0,4){};
      \node[place,tokens=3](B)at(0,2){};
      \node[place,tokens=0](C)at(2,4){};
      \node[place,tokens=2](D)at(2,2){};
      \node[place,tokens=0](G) at(4,3){};
      \node[transition,label=below:$r$](r1)at(1,3){}
      edge[pre]node[above]{$1$}(A)
      edge[post]node[above]{$2$}(C)
      edge[pre]node[above]{$3$}(B);
      \node[transition,label=below:$s$](r3)at(3,3){}
      edge[pre]node[above]{$1$}(C)
      edge[pre]node[above]{$2$}(D)
      edge[post]node[above]{$1$}(G);
    \end{tikzpicture}
  \caption{A Petri net with an initial marking.\label{fig:initial-marking}}
\end{figure}
Over this fixed structure of possible events and conditions, a causal dependency (or flow) relation between sets of events and conditions is described via pre- and post-relations, and it is this structure which determines the possible dynamic behaviour of the net.
A transition in this causal dependency relation is enabled if all places connected to it as inputs contain at least the number of tokens indicated by the arc.
The transition \(r\) in \Cref{fig:initial-marking} is enabled as its preconditions contain enough tokens, while transition \(s\) needs \(r\) to fire first as \(r\) would produce enough tokens for \(s\) to fire.
The dynamic behaviour of a Petri net is often referred to as ``token game''.

\paragraph{Petri nets via Dialectica Categories.}
Petri nets were described categorically in several works~\cite{brown1990,meseguer1990,marti2005petri} and are still been discussed~\cite{pawel2014reachability,baez2020open,master2019}. Models need to capture the practitioner's imagination and make themselves useful, both for calculations and for insights. Categorical models can be useful for both insights and calculations, but we have not seen categorical models that encompass different kinds of transitions in a single net.

Petri nets were modelled using Dialectica categories~\cite{dePaiva1991} previously, but the original Brown and Gurr model~\cite{brown1995} worked only for \textit{elementary nets}, that is nets whose transitions are marked with \{0,1\} for presence or absence of a relationship. An extension of this modelling to deal with integers $\Naturals$ was planned~\cite{dePaiva1991multirelations,brown1991}, but never published. In this work we put together different kinds of transitions, in a single categorical framework. This way the categorical modelling applies to the many kinds of newer applications~\cite{emzivat2016probabilistic} that already use different kinds of labels on the transitions (see \Cref{sec:transition-types}  for a brief account of Petri net transitions often used in applications).

The original dialectica construction~\cite{dePaiva1991} was given in two different styles called the categories $DC$~\cite{dePaiva1989:AMS} and the categories $GC$~\cite{dePaiva1989:CTCS}. For  both constructions $C$ is a  cartesian closed category with some other structure. The first style is connected to G\"odel's Dialectica Interpretation hence the `D' in $DC$ for Dialectica. The second style  called $GC$ (\cite{dePaiva1989:CTCS}) is based on a suggestion of Girard's (hence the `G') on how to simplify the first construction, if one wants a model of Linear Logic.
These two constructions are connected, via monoidal comonads as described in~\cite{dePaiva1991}.
Here we are mostly interested in the construction called $GC$, whose morphisms are simpler. 
When the category $C$ is $\Set$, this construction can be presented in two ways.
This is because a relation in $\Set$ between $U$ and $X$ can be thought of as either a subset of the product, $\alpha \subseteq U\times X$, or as map into \(2\), $\alpha\colon U\times X\to 2$.

The work here uses only the second presentation, defining general relation maps into algebraic structures called \textit{lineales}. This is because changing the lineale where our relations take `values', gives us the possibility of modelling several different kinds of processes. The original dialectica construction deals only with the Heyting algebra-like  lineale $2$. Here we  discuss several other lineales and dialectica categories built over these different lineales.

\eject
We will use some intuitions from the game semantics for linear logic as described by Blass~\cite{blass1992game}, who explicitly compares it with de Paiva's Dialectica interpretation~\cite{dePaiva1991}.

This paper is about models of Petri nets in Dialectica categories.
Petri nets have several different definitions in the literature and multifaceted applications as well.
This work emphasises the underlying graph of a Petri net, which is  a kind of labelled, directed \(L\)-weighted hypergraph, and focuses on a static semantics, simply in terms of ways of putting together combinations of underlying graphs.
This is a first step in understanding these ubiquitous modelling systems, their executions and their dynamic semantics.
While we will suggest semantic interpretations of the different logics for Petri nets, we focus on their compositions via linear logic connectives.
In this, we follow a line of research that explores how to combine Petri nets with linear logic connectives~\cite{brown1995,engberg1990petri}.

\paragraph{Related work.}

Our work fits in the vast landscape of categorical approaches to Petri nets, building on~\cite{brown1991,dePaiva1991multirelations}.
Meseguer and Montanari's seminal work~\cite{meseguer1990} focused on reachability properties of Petri nets, defining a category of all possible executions of a net. This work adopted the collective tokens philosophy. Its ideas were extended to the individual tokens philosophy in~\cite{bruni2001functorial}.
Marti-Oliet and Meseguer explore the connections between linear logic and Petri nets~\cite{marti2005petri}.
Other categorical models of Petri nets focus on obtaining nets by composing smaller nets along some boundaries. One of the first compositional models doing this was~\cite{katis1997nets} where nets are composed along common places. In~\cite{pawel2014reachability}, nets are composed along common transitions and compositionality is used to study reachability properties of Petri nets.
The work of~\cite{brown1990,dePaiva1991multirelations} and~\cite{brown1991} concentrates on combining Petri nets via different monoidal products that give to the category of Petri nets a linear logic structure.
More recently, there has been numerous works building on the ideas of~\cite{meseguer1990} and adopting the formalism of~\cite{katis1997nets}. In~\cite{baez2020open} and~\cite{master2020petri}, the authors focus on studying the categorical properties of reachability. In~\cite{kock2020elements}  a more fine-grained categorical model is proposed, that allows Kock to encompass the individual and collective token philosophies in the same framework. Finally,~\cite{baez2021categories} constructs a unifying framework for~\cite{meseguer1990,bruni2001functorial} and~\cite{kock2020elements} extending~\cite{master2020petri}.

Our work extends the approach of~\cite{brown1990} to allow different kinds of arcs, e.g. inhibitor, probabilistic, partially defined,  natural/integer numbers valued, and the coexistence of them in the same net.

\paragraph{Outline.}
In \Cref{sec:lineales}, we will first provide a brief overview of the types of labels commonly used in Petri net applications. We will use these findings to motivate the need for a categorical model of Petri nets capable of encoding this broad spectrum of relations. We will conclude this section by transforming sets of labels into structures called lineales, whose configuration will later be leveraged to build linear logic connectives.
We show that lineales form a cartesian category.

The construction of our categorical model of Petri nets will take place in two steps. First, in \Cref{sec:intermidiate-cat}, we will build an intermediate category called $\dialLineale{L}$, that will be used to encode pre- and  post-conditions.
Its name is a mnemonic for multisets, where the multitude of elements is specified by values in the lineale \(L\) instead of by just natural numbers.
Within this category, we will define linear connectives and prove that it forms a symmetric monoidal closed category. Subsequently, in \Cref{sec:Petri-net-cat}, we will take two instances of the category $\dialLineale{L}$, one to encode preconditions and the other to encode postconditions, and glue them together (by taking a pullback in the category \(\Cat\)) to build our category \(\dialNet{L}\) of Petri nets. Additionally, in this section, we will develop logical connectives for the category \(\dialNet{L}\) and, in particular, demonstrate that it inherits the symmetric monoidal closed structure from $\dialLineale{L}$.

Finally, \Cref{sec:changing-logic} presents examples of Petri nets spanning a wide spectrum of linear structures, all motivated by various applications, and concludes by showing functoriality of the constructions \(\dialLineale{L}\) and \(\dialNet{L}\).

\section{Lineales: a transition structure for each application of Petri nets}\label{sec:lineales}

In order to define our category of Petri nets, we will take a look at applications to motivate the kind of structure, that of a \textit{lineale}, that will be used to encode a general class of pre- and post-condition relations.
After exposing the motivating application in \Cref{sec:transition-types}, this section presents the definition of lineale (\Cref{sec:lineales-def}) and gives examples of them (\Cref{sec:lineales-examples}).

\subsection{Chemical and metabolic networks: an inspiration}\label{sec:transition-types}

Networked systems are determined by their connections~\cite{Jost2020}.  Perhaps the most basic type of relationship in any network is one that only allows us to express either presence or absence, that is, where the relationship connecting nodes uses the set \{0,1\} as a ruler or label set.  In real-world applications this is, though, not sufficient.  In this section we explore frequent and rich applications of Petri nets, from chemical reaction networks to metabolic networks, searching for the kind of labels used on pre- and post- conditions.

\textit{Chemical reaction networks.} Chemical combination is compositional in nature.  Although data on substance reactivity are typically annotated as a list of chemical equations (\Cref{ChemDataRep}, left), chemists reason on the network structure (\Cref{ChemDataRep}, right) that emerges when the reactions are connected to make their concurrency explicit~\cite{Schummer1998}.  Synthesis planning is a prominent example: suppose we want to synthesise substance F, but we cannot carry out reaction $r_2$, because we have no substance C  on our lab's shelf. In such a case, equation $r_1$ provides another synthetic route, since it is possible to obtain F from A and B (and E). In other words, the synthesis of F results from composing reactions $r_1$ and $r_2$.

Directed hypergraphs and their enhancements, such as Petri nets, are used to model chemical reaction networks for they are models of concurrency of directed relations.  These models provide a rich semantic basis on which to interpret questions that arise in chemistry, such as what substances can be synthesised from a given set of starting materials?~\cite{Stadler2018}
Given a target substance, which synthetic routes are known and which starting materials are needed to reach the target?~\cite{Hoffmann2009} Do chemical reactions turn targets into key precursors?~\cite{Jost2020}
How many synthetic routes pass through a given reaction? These questions can be answered by probing the topology and geometry of the wiring of the network. The first two questions are answered by defining suitable closure operators~\cite{Stadler2018}, and the last two questions by computing the curvature of the edges of the network~\cite{forman2018}, taking into consideration the proportions in which substances combine (stoichiometric coefficients).

At the level of abstraction described above, transitions of chemical reaction networks are discrete in nature, and pre- and post- conditions correspond to presence/absence of substances or to stoichiometric coefficients, which can be modelled by the rulers $\{0,1\}$ and $\Naturals$, respectively.

\textit{Metabolic networks.} These networks comprise the metabolic pathways (network of chemical reactions) and the gene interactions that regulate them.  A key aspect of the former is the kinetic modelling.  There, Petri nets model reaction rates.  For elementary reactions, which take place in a single step, the Law of Mass Action states that reaction rates are proportional to the concentration of reactants.  Both quantities, rate of reactions and concentration of reactants, are usually taken as positive real numbers; therefore, in this application, Petri nets are challenged to handle continuous tokens and transitions, which requires the ruler $\Reals^+$. On the other hand, gene interactions are handled by implementing genetic switches that are modelled by discrete transitions.  A Petri net model for a metabolic network therefore needs two different rulers on the same net.

When applied to concrete metabolisms, a Petri net model will usually need to incorporate more than two rulers at the same time.  For instance,~\cite{Saito2011} shows a hybrid Petri net representation of the gene regulatory network of \textit{C. elegans} that is labelled with discrete and continuous transitions, but also with  negative integers, real numbers, strings, and products of them.

Summarising,  applications may need rulers such as $\{0,1\}$, $\{-1,0,-1\}$ (for data uncertainty, which is common in complex network systems~\cite{PhysRevResearch.2.043135}), $\Naturals$, $\Reals^+$, $\Integers, \Reals$, strings,  and their finite products.  The ability of choosing from a vast pool of rulers to label pre- and post- conditions is one of the strengths of the categorical construction presented in this paper.

\subsection{Lineales: the codomain of general relations}\label{sec:lineales-def}

% \elena{why do we want lineales?
% which part of the classical story we take as fundamental?
% disjunction and the weird coimplication?
% the interpretation seems weird}
Classical relations are functions of type \(U \times X \to 2\).
If a relation assigns \(1\) to a pair \((u,x)\), then the two elements are related, otherwise they are not.
In some applications, one would like more information about the relationship between the pair \((u,x)\).
For example, we could record the intensity of the relation by assigning to each pair \((u,x)\) a natural number.
In general, we can define poset-valued relations.
For our purposes, posets will have additional structure: a multiplication and an internal-hom.
The multiplication makes the poset also a monoid.
\Cref{sec:intermidiate-cat} uses this structure to define a monoidal closed structure on lineale-valued relations (\Cref{th:monclosed}).

\begin{definition}[Partially ordered monoid]
  A \emph{partially ordered commutative monoid} $(L,\linorder,\lintensor,e)$ is a commutative monoid $(L,\lintensor, e)$ equipped with a partial order $\linorder$ that is compatible with the monoidal operation, i.e.~if $a \linorder b$ and \(a' \linorder b'\) then $a \lintensor a' \linorder b \lintensor b'$.
\end{definition}
The internal-hom is the best approximation to an inverse operation for the multiplication of the monoid.
\begin{definition}[Internal-hom in a monoid]
  Let $(L,\linorder,\lintensor,e)$ be a partially ordered commutative monoid. A binary operation \(\lininthom \colon L^{op} \times L \to L\) is said to be an \emph{internal-hom} when it is right adjoint to the monoidal product \(\lintensor\), i.e. \(\forall a,b,c \in L, \ b \lintensor c \linorder a \Leftrightarrow b \linorder c \lininthom a\).
  The internal-hom is also required to respect the ordering, contravariantly in the first coordinate and covariantly in the second, i.e.~if \(b \linorder a\) and \(a' \linorder b'\) then \(a \lininthom a' \linorder b \lininthom b'\).
\end{definition}
A lineale~\cite{dePaiva2002} is a monoidal closed poset.
This means that a lineale has a commutative monoid structure with a right adjoint, the internal-hom, that are both compatible with the order.  
This structure is also known as a residuated ordered commutative monoid~\cite{Hohle1995}.

\begin{definition}
  A \emph{lineale} is a tuple $(L,\linorder,\lintensor,e,\lininthom)$ such that $(L,\linorder,\lintensor,e)$ is a partially ordered monoid and $\lininthom$ is an internal-hom for $(L,\linorder,\lintensor,e)$.
\end{definition}
% \begin{example}\label{po-monoid_lineales}
%   Examples of lineales are $(\Naturals,\geq,+,0,\lininthom)$ and $(\Reals^+,\geq,+,0,\lininthom)$,
%   where\\
%   \(a \lininthom b=\begin{cases} b - a &  a\leq b\\ 0 &  a>b \end{cases}\).
% \end{example}

Notice that in any lineale, \(b = e \lininthom b\) for any \(b \in L\). Since \(b \lintensor e \linorder b\) holds in any monoid \((L, \lintensor, e)\), it follows by the definition of \(\lininthom\) that \(b \linorder e \lininthom b\) for all \(b \in L\).   Moreover, choosing \(c = e \lininthom b\) in $c \lintensor e \linorder b \iff c \linorder e \lininthom b$ we obtain $e \lininthom b = (e \lininthom b) \lintensor e \linorder b \iff e \lininthom b \linorder e \lininthom b$. Since the right side of the equivalence is trivially true, we must have  $e \lininthom b \linorder  b$.

\begin{definition}
  A \emph{morphism of lineales} \(h \colon (L,\linorder,\lintensor,e,\lininthom) \to (L',\linorder',\lintensor',e',\lininthom')\) is a monotone function \(h \colon (L,\linorder) \to (L',\linorder')\) that laxly preserves the monoid structure, i.e.\ \(e' \linorder' h(e)\) and \(h(a) \lintensor' h(b) \linorder' h(a \lintensor b)\).
\end{definition}

In the same way that lineales are posetal monoidal closed categories, morphisms of lineales are lax monoidal functors between them.
Lineales and their morphisms form a category \(\Lin\).
Note that the preservation of the internal-hom follows from the adjointness property:
since \(h(a \lininthom b) \lintensor' h(a) \linorder' h((a \lininthom b) \lintensor a) \linorder' h(b)\), then \(h(a \lininthom b) \linorder' h(a) \lininthom' h(b)\).

The next lemma recalls the relationship between minima and internal-homs.
It is a poset instance of the well-known adjoint functor theorem, see, for instance, Section V.6 in~\cite{maclane}.

\begin{lemma}\label{lemma:finite-pomonoids}
  A partially ordered monoid \((M,\linorder,\lintensor,e)\) in which infima always exist defines a lineale with internal-hom \(a \lininthom b = \inf\{x \in M \mid x \lintensor a \linorder b\}\).
\end{lemma}
% \begin{proof}
%   Suppose that \(a \lintensor b \linorder c\).
%   Then, by definition of \(\lininthom\), \(a \linorder b \lininthom c\).
%   For the converse, suppose that \(a \linorder b \lininthom c\).
%   By definition of \(\lininthom\) and existence of minima, \((b \lininthom c) \lintensor b \linorder c\).
%   Then, \(a \lintensor b \linorder (b \lininthom c) \lintensor b \linorder c\) because \(\lintensor\) is compatible with \(\linorder\).

%   We check that \(\lininthom\) is compatible with the ordering.
%   Suppose that \(b \linorder a\) and \(a' \linorder b'\).
%   By the definition of \(\lininthom\), \((a \lininthom a') \lintensor a \linorder a'\).
%   Then, \((a \lininthom a') \lintensor b \linorder (a \lininthom a') \lintensor a \linorder a' \linorder b'\) because the multiplication is compatible with the ordering and by the assumption.
%   Then, \(a \lininthom a' \linorder b \lininthom b'\) because of the definition of internal-hom.
% \end{proof}
For partially ordered groups, the internal-hom is exactly the inverse operation to the multiplication.
In fact, any partially ordered group is a lineale with $a\lininthom b= b \lintensor a^{-1}$ (see also~\cite[Example 4.4.1]{Seal_Tholen_2014}).
Some of our examples will fall into this case.

\begin{definition}[Partially ordered group]\label{monoidal_poset}
  A \emph{partially ordered group} $(G,\linorder,\lintensor,e,(-)^{-1})$ is a partially ordered monoid \((G,\linorder,\lintensor,e)\) together with an inverse operation \((-)^{-1}\) that makes \((G,\lintensor,e,(-)^{-1})\) a group and respects the ordering contravariantly, i.e.~if \(a \linorder b\) then \(b^{-1} \linorder a^{-1}\).
\end{definition}

\begin{lemma}\label{po-group_lineale}
  A partially ordered group $(G,\linorder,\lintensor,e,(-)^{-1})$ can be endowed with the structure of a lineale with $a \lininthom b \coloneqq b \lintensor a^{-1}$.
\end{lemma}
\begin{shortproof}
  A group is a monoid, thus we only need to check that \(a \lininthom b\) actually defines an internal-hom and that it respects the ordering.
  \begin{align*}
    & b \lintensor c \linorder a & & b \linorder a \land a' \linorder b'\\
    & \Leftrightarrow b \lintensor c \lintensor c^{-1} \linorder a \lintensor c^{-1} & & \Rightarrow a^{-1} \linorder b^{-1} \land a' \linorder b'\\
    & \Leftrightarrow b \linorder c \lininthom a & &\Rightarrow a' \lintensor a^{-1} \linorder b' \lintensor b^{-1}\\
    && & \Rightarrow a \lininthom a' \linorder b \lininthom b'
  \end{align*}
\end{shortproof}

% \begin{example}
%   Examples of lineales obtained from partially ordered groups are $(\Integers,\geq,+,0,-)$ and $(\Reals,\geq,+,0,-)$.
% \end{example}

\subsection{Examples of lineales}\label{sec:lineales-examples}
While the lineale $2$ is associated with Boolean and Heyting algebras, which are traditional algebraic models for classical and intuitionistic propositional logic, other lineales are associated with different non-classical systems.
We describe some lineales and the non-classical logic associated to them.
Finally, we consider the coexistence of several lineales in a single net by taking products of them.

\begin{example}[Classical lineale]\label{ex:classical-lineale}
  The original work on the categorical version of the Dialectica interpretation has concentrated on relations that take values into $2$, considered as a lineale.
  The set \(2 = \{0,1\}\) has a lineale structure \((2,\geq,\lor,0,\lininthom)\), where \(1 \geq 0\), the operation \(\lor\) is the logical disjunction if \(0\) is interpreted as \emph{false} and \(1\) as \emph{true}, and the internal-hom is defined as \(a \lininthom b = \min\{x \in 2 \mid x \lor a \geq b\}\).
  This is a lineale structure by \Cref{lemma:finite-pomonoids}.
  More explicitly, \(a \lininthom b = 0\) if \(a \geq b\) and \(a \lininthom b = b\) otherwise.
  This operation is the adjoint of the disjunction and differs from the logical implication, which is adjoint to the conjunction.
  In fact, \(a \lininthom b = \lnot (b \to a)\).
\end{example}

\begin{example}[Kleene lineale]\label{lem:kleene-lineale}
  We describe a 3-valued propositional logic where the undefined truth-value, the  ``unknown'' state, can be thought of as neither true nor false.
  We interpret \(0\) and \(1\) as \emph{false} and \emph{true} respectively.
  The additional truth value \(u\) represents \emph{undefined}.
  The three elements set \(3=\{0,u,1\}\) has a lineale structure \((3,\geq,\max,0,\lininthom)\), where \(1 \geq u \geq 0\) and the internal-hom is defined as \(a \lininthom b = \min\{x \in 3 \mid \max\{x,a\} \geq b\}\).
  More explicitly, \(a \lininthom b = 0\) if \(a \geq b\) and \(a \lininthom b = b\) otherwise.
  This is indeed a lineale by \Cref{lemma:finite-pomonoids}.
\end{example}
% \elena{why do we mention this? should we make this into an example then?}
We should  also mention the lineale $4$, associated with Belnap-Dunn's four-valued logic~\cite{belnap1977useful}.
These four values also correspond to the algebraic identities for the two conjunctions and two disjunctions of Linear Logic.

\begin{example}[Multirelations lineale]\label{lem:naturals-lineale}
  We consider a lineale defined on the natural numbers, \((\Naturals,\geq,+,0,\lininthom)\), where the order is the usual order on natural numbers, \(m \linorder n\) iff \(m \geq n\), and the monoid structure is that of addition.
  While we can think of the classical truth values as indicating whether or not two elements are related, we can think of the natural numbered truth values as indicating how many times two elements are related.
  Note that this lineale is not a quantale as suprema need not exist.
  % This idea is related to Lawvere's generalised metric spaces~\cite{lawvere1973}.
  The internal-hom is defined as in \Cref{lemma:finite-pomonoids} and gives a lineale structure.
  In the case of natural numbers, this internal-hom, \(a \lininthom b = \min\{n \in \Naturals \mid n + a \geq b\}\), becomes \(b-a\), when \(b \geq a\), and \(0\), otherwise.
  In other words, \(a \lininthom b = \max \{0, b-a\}\).
\end{example}

\begin{example}[Integers lineale]\label{ex:integers-lineale}
  % We believe this style of dialectica category can be profitably used to model systems where some transitions can cancel other transitions.
  % \elena{what is the interpretation of integers-valued relations? what does it mean to be negatively related?}
  Similarly to the multirelations lineale, we consider a lineale structure on the integers, \((\Integers, \geq, +, 0, \lininthom)\).
  As the monoid of integers with addition is actually a group, we can apply \Cref{po-group_lineale} to choose the internal-hom structure as subtraction: \(a \lininthom b = b - a\).
\end{example}

\begin{example}[Probabilistic lineale]\label{lem:prob-lineale}
Next we consider the reals, in the form of the closed interval $[0,1]$. These have long been considered for fuzzy sets, as the real number associated with a pair $(u,x)$ can be thought of as the probability of the association between $u$ and $x$.
% So far we considered only (topologically) discrete lineales. Now we want to discuss real numbers, in particular the closed interval $[0,1]$.
We show that the closed interval $[0,1]$ admits a lineale structure.
\begin{itemize}
    \item The monoid structure is given by the product of real numbers, \(a \lintensor b \coloneqq a \cdot b\), whose unit is \(1\);
    \item The partial order is the usual ordering on the reals;
    \item The internal-hom is a `truncated division' \(a \lininthom b \coloneqq \frac{b}{a}\) if \(a > b\) and \(a \lininthom b \coloneqq 1\) otherwise.
\end{itemize}
This structure defines a lineale by \Cref{lemma:finite-pomonoids}: if \(a > b\), then \(\inf\{x \in [0,1] \mid x \cdot a \geq b\} = \frac{b}{a}\), while if \(a \leq b\), then \(\inf \emptyset = 1\).
\end{example}
% \begin{lemma}
% The closed interval \([0,1]\) is a lineale with the structure defined above.
% \end{lemma}
% \begin{proof}
% The image of the product of real numbers in the closed interval \([0,1]\) is the closed interval itself. This product is associative and unital, thus it gives a monoid structure to \([0,1]\).
% Moreover, it preserves the ordering and this makes \([0,1]\) a partially ordered monoid.
% We need to check that truncated division as defined above gives an internal-hom.
% We reason by cases as our metatheory is classical.
% \begin{align*}
%     &\text{If }c \neq 0 \land c \geq a &&\text{If }c = 0 &&\text{If }c<a\\
%     & b \linorder c \lininthom a && b \linorder c \lininthom a && b \linorder c \lininthom a \\
%     \Leftrightarrow \ &b \leq \dfrac{a}{c} &\Leftrightarrow \ & b \leq 1 &\Leftrightarrow \ & b \leq 1\\
%     \Leftrightarrow \ &b \cdot c \leq a &\Leftrightarrow \ & 0 \leq a &\Leftrightarrow \ & b \leq 1 < \dfrac{a}{c}\\
%     \Leftrightarrow \ &b \lintensor c \linorder a &\Leftrightarrow \ & b \cdot 0 \leq a &\Leftrightarrow \ & b \cdot c \leq a\\
%     & &\Leftrightarrow \ & b \lintensor c \linorder a &\Leftrightarrow \ & b \lintensor c \linorder a\\
% \end{align*}
% \end{proof}

\begin{example}[Morphism of lineales]\label{ex:lineale-morph}
  Consider the Multirelations lineale from \Cref{lem:naturals-lineale} and the Classical lineale from \Cref{ex:classical-lineale}.
  The function \(h \colon \Naturals \to 2\), defined by \(h(0) = 0\) and \(h(n) = 1\) for all \(n \geq 1\), is a morphism of lineales.
  In fact, it is monotone, if \(m \leq n\) then \(h(m) \leq h(n)\), and it preserves the unit, \(h(0) = 0\), and the multiplication, \(h(m+n) = h(m) \lor h(n)\).
  If we interpret the logical values in \(\Naturals\) as recording `how many times' a proposition is true, this morphism forgets this number and only records whether or not the proposition is true.
\end{example}

\paragraph{Product of lineales.}
% \elena{we could probably show that the category of lineales is complete and cocomplete (or just products and coproducts). so we can take the coproduct of two natural-numbered lineales to get different kinds of arcs}
We have produced a pool of lineales, each of them suitable for  transitions taking values in certain data types ($\{0,1\}$, $\{-1,0,-1\}$, $\Naturals$, $\Integers$ or $\Reals^+$).
As discussed in \Cref{sec:lineales}, in empirical data analysis, a transition often carries data on more than one variable simultaneously.
In this section, we show that any finite combination of lineales can be endowed with the structure of a lineale by taking finite products of them.

As lineales are just the poset-version of symmetric monoidal closed categories, the next proposition is the posetal version of the analogous result for symmetric monoidal closed categories.
Symmetric monoidal closed categories form a category $\textbf{SymClosedCat}$ whose objects are symmetric monoidal closed categories, and morphisms are functors that preserve the adjunction.
It is folklore that this category is cartesian.
For completeness, we directly show the particular case of our interest: the component-wise product of two lineales is a lineale.

\begin{proposition}\label{prop:product-lineale}
  If $(L_1,\leq_1,\lintensor_1,e_1,\lininthom_1)$ and $(L_2,\leq_2,\lintensor_2,e_2,\lininthom_2)$ are lineales, then \(L_{1} \times L_{2}\) has a lineale structure where, for elements $l=(l_1,l_2), l'=(l'_1,l'_2) \in L_1\times L_2$,
  % $(L_1\times L_2,\leq,\lintensor,e,\lininthom)$ is a lineale with the following structure.
  \begin{itemize}
    \item the product is defined component-wise, $l\lintensor l' = (l_1\lintensor_1 l'_1, l_2\lintensor_2 l'_2)$;
    \item the unit is the pair of units, $e=(e_1,e_2)$;
    \item the order $l\leq l'$ holds if and only if $l_1\leq_1 l'_1$ and $l_2\leq_2 l'_2$; and
    \item the internal-hom is defined component-wise, $l \lininthom l'= (l_1\lininthom_1 l_1',l_2\lininthom_2 l_2')$.
  \end{itemize}
  This operation is the product in the category of lineales, \(\Lin\).
\end{proposition}
\begin{proof}
  $(L_1\times L_2,\lintensor,e)$ is the cartesian product of two monoids and therefore it is a monoid. $(L_1\times L_2,\leq)$ is a partial ordered set with the ordering defined above.
  Since $l_i \leq l'_i$ implies both $l_i\lintensor_{i} k_i \leq l'_i\lintensor_{i} k_i$ and $k_i\lintensor_{i} l_i \leq k_i\lintensor_{i} l'_i$ for each $l_{i}, k_i\in L_i$ and $i = {1,2}$; then $l\leq l'$ implies $l\lintensor k \leq l'\lintensor k$ and $k\lintensor l \leq k\lintensor l'$ for every $k=(k_1,k_2) \in L_1\times L_2$. This proves that $(L_1\times L_2,\leq,\lintensor,e,\lininthom)$ is a partially ordered monoid.
  We need to prove that the internal-hom defined above is right adjoint to the monoidal product.
  \begin{align*}
    & b \lintensor c \leq a\\
    \Leftrightarrow \ & (b_1 \lintensor_1 c_1, b_2 \lintensor c_2) \leq (a_1,a_2)\\
    \Leftrightarrow \ & b_1 \lintensor_1 c_1 \leq_1 a_1 \land b_2 \lintensor_2 c_2 \leq_2 a_2\\
    \Leftrightarrow \ & b_1 \leq_1 c_1 \lininthom_1 a_1 \land b_2 \leq_2 c_2 \lininthom_2 a_2\\
    \Leftrightarrow \ & (b_1,b_2) \leq (c_1 \lininthom_1 a_1,c_2 \lininthom_2 a_2) \\
    \Leftrightarrow \ & b \leq c \lininthom a
  \end{align*}
  This proves that the product of lineales is again a lineale.

  Finally, we check the universal property.
  The projections of the underlying sets are also strict morphisms of lineales, as the structure is defined component-wise.
  For two morphisms of lineales \(g \colon (L, \linorder, \lintensor, e, \lininthom) \to (L_{1}, \linorder_{1}, \lintensor_{1}, e_{1}, \lininthom_{1})\) and \(h \colon (L, \linorder, \lintensor, e, \lininthom) \to (L_{2}, \linorder_{2}, \lintensor_{2}, e_{2}, \lininthom_{2})\), there is a unique function \(\langle g, h \rangle \colon L \to L_{1} \times L_{2}\) commuting with the projections by the universal property of the product in \(\Set\).
  This function also preserves the order and the monoidal structure as it does so component-wise.
\end{proof}

%%%

\section{The category $\dialLineale{L}$}\label{sec:intermidiate-cat}
Having defined a lineale, we proceed to construct the intermediate category $\dialLineale{L}$ over which our category of Petri nets \(\dialNet{L}\) is built.

Its objects \((U,X,\alpha)\) represent \(L\)-valued relations \(\alpha\) between the sets \(U\) and \(X\).
The game semantics interpretation of \(\alpha\) is given by the game played on the proposition \(\bigvee_{u \in U} \bigwedge_{x \in X} \alpha(u,x)\): Proponent tries to prove \(\alpha\) and chooses an element \(u \in U\), then Opponent tries to refute it and chooses an element \(x \in X\); the result of the game is determined by the value of \(\alpha(u,x)\).
In the classical case, \(L = \{0,1\}\), Proponent wins if \(\alpha\) is true on the chosen pair \((u,x)\), \(\alpha(u,x) = 1\).
In the other cases, we may imagine that the possible outcomes of the game are not just ``win'' and ``loose'', but they carry the extra information given by the lineale.

Morphisms need to carry the information about both the players moves~\cite{blass1992game}.
Thus, morphisms \((U,X,\alpha) \to (V,Y,\beta)\) are pairs of functions, the first one \(f \colon U \to V\) that maps every choice of Proponent in \(\alpha\) to one in \(\beta\) and the second one \(F \colon Y \to X\) that maps every choice of Opponent in \(\beta\) to one in \(\alpha\).
Suppose that Proponent has a winning strategy in \(\alpha\), i.e.~there is \(u \in U\) such that, for all \(x \in X\), \(\alpha(u,x)\) is true.
Then, Proponent can choose \(f(u)\) in \(\beta\) and we know that each choice \(y \in Y\) of Opponent will lead to an outcome \(\beta(f(u),y)\) that is bounded by \(\alpha(u,F(y))\).

These morphisms make \(\dialLineale{L}\) a category (\Cref{prop:category}).

\begin{definition}[Category \(\dialLineale{L}\)]
Given a lineale $(L,\linorder,\lintensor,e,\lininthom)$, the category $\dialLineale{L}$ is defined by the following data.
\begin{itemize}
    \item An object is a triple $A = (U,X, \alpha)$, denoted by $U\dialarrow{\alpha}X$, where $U, X$ are sets and $\alpha \colon U \times X \to L$ is a function in $\Set$.
    \item  A morphism \((f,F) \colon (U \dialarrow{\alpha} X) \to (V \dialarrow{\beta} Y)\) is a pair of morphisms, $f\colon U \to V$ and $F\colon Y \to X$ in $\Set$, such that $\forall u \in U \ \forall y \in Y \ \alpha(u, Fy) \linorder \beta(fu,y)$.
\end{itemize}
\begin{equation*}
  \begin{tikzcd}[]
    & U \times Y \arrow[rr,"f \times \id{Y}"] \arrow[dd,"\id{U} \times F"']& & V \times Y \arrow[dd,"\beta"] \\
    & & \linorder & \\
    & U \times X \arrow[rr,"\alpha"'] & & L \\
  \end{tikzcd}
\end{equation*}
\end{definition}

The category $\dialLineale{L}$ allows us to have $L$-valued relations, including multirelations ($L=\Naturals$) and any other label set that can be seen as a lineale.
The objects of \(\dialLineale{L}\) can be seen as \(L\)-enriched profunctors~\cite[Section 7.7]{Borceux_1994}, also known as distributors, where the lineale \(L\) is seen as a one-object monoidal closed category.
Their morphisms are, usually, natural transformations, while here we consider lax ones.

\begin{example}\label{ex:dialectica-objects-morphisms}
  For three elements \(a,a',b \in L\) of a lineale \((L,\linorder,\lintensor,e,\lininthom)\) and sets \(U = \{u,u'\}\), \(X = \{x,x'\}\), \(V = \{v\}\) and \(Y = \{y\}\), we can define two objects \((U \dialarrow{\alpha} X)\) and \((V \dialarrow{\beta} Y)\) of \(\dialLineale{L}\).
  \begin{align*}
    \alpha(u,x) & = a & \alpha(u',x) & = a'\\
    \alpha(u,x') & = e & \alpha(u',x') & = e\\
    \beta(v,y) & = b & &
  \end{align*}
  These can be represented as weighted bipartite graphs where the elements of \(U\) and \(X\) are the vertices and are connected by an edge with weight \(a \in L\) whenever the value of \(\alpha\) on them is \(a\) (\Cref{fig:weighted-relations-morphisms}).
  Omitted edges implicitly have weight \(e\), the unit of the lineale \(L\).
  \begin{figure}[h!]
    \centering
    \linealeRelationMorphismFig{}
    \caption{Two morphisms in \(\dialLineale{L}\), if \(a,a' \linorder b\) (left) and if \(b \linorder a,e\) (right).\label{fig:weighted-relations-morphisms}}
  \end{figure}
  If \(a,a' \linorder b\), there is a morphism \((f,F) \colon (U \dialarrow{\alpha} X) \to (V \dialarrow{\beta} Y)\) and, if \(b \linorder a,e\), there is a morphism \((g,G) \colon (V \dialarrow{\beta} Y) \to (U \dialarrow{\alpha} X)\).
  \begin{align*}
    f(u) = f(u') & = v & F(y) & = x\\
    g(v) & = u & G(x) = G(x') & = y
  \end{align*}
\end{example}

\begin{proposition}\label{prop:category}
  $\dialLineale{L}$ is a category.
\end{proposition}
\begin{proof}
  The identity arrow of an object $U\dialarrow{\alpha}X$ in $\dialLineale{L}$ is given by the pair $(\id{U},\id{X})$ of identities in $\Set$.  Moreover, given objects $A=(U\dialarrow{\alpha}X)$,  $B=(V\dialarrow{\beta}Y)$, and $C=(W\dialarrow{\gamma}Z$), and morphisms $(f,F)\colon A\to B$ and $(g,G)\colon B \to C$, their composition is computed componentwise as $(g,G)\compose (f,F)=(g\compose f,F\compose G)\colon A\to C$. Notice that $(g,G)\compose (f,F)$ is a morphism in $\dialLineale{L}$: given $u \in U$ and $z\in Z$, we have $\alpha(u,FGz)\linorder \beta(fu,Gz)\linorder \gamma(gfu,z)$. Associativity and unitality come from associativity and unitality in $\Set$.
\end{proof}

We now proceed to define products and coproducts in \(\dialLineale{L}\) and to equip it with a symmetric monoidal closed structure. To achieve this, we lift the symmetric monoidal closed structure of the lineale: the product and coproduct rely only on the set structure, the tensor product is defined using the monoidal structure of the lineale, and the internal-hom in \(\dialLineale{L}\) is derived from the internal-hom of the lineale.

% We will define maps up to symmetries in \(\Set\) to avoid distracting the reader with details.
The intuitions from the game semantics for linear logic justify the types of the product and coproduct in the Dialectica construction~\cite{blass1992game}.
For Proponent to win a game \(A \product B\), it needs to have a winning strategy both in \(A\) and in \(B\), while Opponent only needs to choose one between \(A\) and \(B\), and have a winning strategy there.
Thus, in the Proponent side, the product of \(A\) and \(B\) is the categorical product \(U \times V\), while in the Opponent side, it is the coproduct \(X + Y\).
Similarly, for Proponent to win a game \(A \coproduct B\), it needs to choose one between \(A\) and \(B\), and have a winning strategy there, while Opponent needs to have a winning strategy both in \(A\) and in \(B\).
This is reflected in the types of the coproduct, \(U + V\) for the Proponent side and \(X \times Y\) for the Opponent side.

\begin{definition}[Product and coproduct in \(\dialLineale{L}\)]\label{def:cartesian}
  Given two objects $A = (U \dialarrow{\alpha} X)$ and $B = (V \dialarrow{\beta} Y)$ in $\dialLineale{L}$, we define their cartesian product $A \product B$ and their coproduct \(A \coproduct B\) as the following objects.
  \begin{align*}
    A \product B &= (U \times V \dialarrow{\alpha \product \beta} X + Y) & A \coproduct B = (U + V \dialarrow{\alpha \coproduct \beta} X \times Y)
  \end{align*}
  The function \(\alpha \product \beta\) is
  \(U \times V \times (X+Y) \xrightarrow{\coproductmap{\alpha \compose \discard_V}{\beta \compose \discard_U}} L \),
  where \(\discard_U \colon U \times V \times Y \to V \times Y\) is the function that discards \(U\) and \(\discard_{V} \colon U \times V \times X \to U \times X\) the function that discards \(V\) in \(\Set\).
  Similarly, the function \(\alpha \coproduct \beta\) is
  \((U + V) \times X \times Y \xrightarrow{\coproductmap{\alpha \compose \discard_Y}{\beta \compose \discard_X}} L \),
  where \(\discard_Y \colon U \times X \times Y \to U \times X\) is the function that discards \(Y\) and \(\discard_{X} \colon V \times X \times Y \to V \times Y\) the function that discards \(X\) in \(\Set\).
\end{definition}

\begin{proposition}\label{prop:dialectica-co-products}
  The operations in \Cref{def:cartesian} give products and coproducts in \(\dialLineale{L}\).
\end{proposition}
\begin{proof}
  The projections \(\pi_{A} \colon A \product B \to A\) and \(\pi_{B} \colon A \product B \to B\) are defined by projections and coprojections in \(\Set\): \(\pi_{A} = (\pi_{U}, \iota_{X})\) and \(\pi_{B} = (\pi_{V}, \iota_{Y})\).

  Let \((f,F) \colon C \to A\) and \((g,G) \colon C \to B\) be two morphisms in \(\dialLineale{L}\).
  We check that the pair \((\productmap{f}{g}, \coproductmap{F}{G})\) is a morphism of type \(C \to A \product B\), using the properties of products and coproducts in \(\Set\) and the definition of morphism in \(\dialLineale{L}\).
  \begin{align*}
    & (\alpha \product \beta) \compose (\productmap{f}{g} \times \id{X+Y})\\
    & = \coproductmap{\alpha}{\beta} \compose ((\productmap{f}{\discard_{V} \compose g} \times \id{X}) + (\productmap{\discard_{U} \compose f}{g} \times \id{Y}))\\
    & = \coproductmap{\alpha}{\beta} \compose ((f \times \id{X}) + (g \times \id{Y}))\\
    & = \coproductmap{\alpha \compose (f \times \id{X})}{\beta \compose (g \times \id{Y})}\\
    & \sqsubseteq \coproductmap{\gamma \compose (\id{W} \times F)}{\gamma \compose (\id{W} \times G)}\\
    & = \gamma \compose (\id{W} \times \coproductmap{F}{G})
  \end{align*}
  The projections composed with \((\productmap{f}{g}, \coproductmap{F}{G})\) give back \((f,F)\) and \((g,G)\).
  \begin{align*}
    & \pi_{A} \compose (\productmap{f}{g}, \coproductmap{F}{G}) && \pi_{B} \compose (\productmap{f}{g}, \coproductmap{F}{G})\\
    & = (\pi_{U} \compose \productmap{f}{g}, \coproductmap{F}{G} \compose \iota_{X}) && = (\pi_{V} \compose \productmap{f}{g}, \coproductmap{F}{G} \compose \iota_{Y})\\
    & = (f,F) && = (g,G)
  \end{align*}
  This is the unique morphism doing so: any morphism \((h,H) \colon C \to A \product B\) that commutes with the projections needs to be \((\productmap{f}{g}, \coproductmap{F}{G})\) by the universal property of the product and coproduct in \(\Set\).
  This shows that \((\product)\) is the categorical product in \(\dialLineale{L}\).

  Dually, we show that \((\coproduct)\) is the categorical coproduct.
  The projections \(\iota_{A} \colon A \to A \coproduct B\) and \(\iota_{B} \colon B \to A \coproduct B\) are defined by projections and coprojections in \(\Set\): \(\iota_{A} = (\iota_{U}, \pi_{X})\) and \(\iota_{B} = (\iota_{V}, \pi_{Y})\).

  Let \((f,F) \colon A \to C\) and \((g,G) \colon B \to C\) be two morphisms in \(\dialLineale{L}\).
  We check that the pair \((\coproductmap{f}{g}, \productmap{F}{G})\) is a morphism of type \(A \coproduct B \to C\), using the properties of products and coproducts in \(\Set\) and the definition of morphism in \(\dialLineale{L}\).
  \begin{align*}
    & (\alpha \coproduct \beta) \compose (\id{U + V} \times \productmap{F}{G})\\
    & = \coproductmap{\alpha}{\beta} \compose ((\id{U} \times \productmap{F}{\discard_{Y} \compose G}) + (\id{V} \times \productmap{\discard_{X} \compose F}{G}))\\
    & = \coproductmap{\alpha}{\beta} \compose ((\id{U} \times F) + (\id{V} \times G))\\
    & = \coproductmap{\alpha \compose (\id{U} \times F)}{\beta \compose (\id{V} \times G)}\\
    & \linorder \coproductmap{\gamma \compose (f \times \id{Z})}{\gamma \compose (g \times \id{Z})}\\
    & = \gamma \compose (\coproductmap{f}{g} \times \id{Z})
  \end{align*}
  The coprojections composed with \((\coproductmap{f}{g}, \productmap{F}{G})\) give back \((f,F)\) and \((g,G)\).
  \begin{align*}
    & (\coproductmap{f}{g}, \productmap{F}{G}) \compose \iota_{A} && (\coproductmap{f}{g}, \productmap{F}{G}) \compose \iota_{B}\\
    & = (\coproductmap{f}{g} \compose \iota_{U}, \pi_{X} \compose \productmap{F}{G}) && = (\coproductmap{f}{g} \compose \iota_{V}, \pi_{Y} \compose \productmap{F}{G})\\
    & = (f,F) && = (g,G)
  \end{align*}
  This is the unique morphism doing so: any morphism \((h,H) \colon A \coproduct B \to C\) that commutes with the coprojections needs to be \((\coproductmap{f}{g}, \productmap{F}{G})\) by the universal property of the product and coproduct in \(\Set\).
\end{proof}

\begin{example}\label{ex:dialectica-objects-co-product}
  The product and coproduct of the objects $A = (U \dialarrow{\alpha} X)$ and $B = (V \dialarrow{\beta} Y)$ defined in \Cref{ex:dialectica-objects-morphisms} are below and in \Cref{fig:dialectica-objects-co-product}.
  \begin{align*}
    \alpha \product \beta (u,v,x) & = a & \alpha \coproduct \beta (u,x,y) & = a \\
    \alpha \product \beta (u',v,x) & = a' & \alpha \coproduct \beta (u,x',y) & = e \\
    \alpha \product \beta (u,v,x') & = e & \alpha \coproduct \beta (u',x,y) & = a' \\
    \alpha \product \beta (u',v,x') & = e & \alpha \coproduct \beta (u',x',y) & = e \\
    \alpha \product \beta (u,v,y) & = b  & \alpha \coproduct \beta (v,x,y) & = b \\
    \alpha \product \beta (u',v,y) & = b & \alpha \coproduct \beta (v,x',y) & = b
  \end{align*}
  \begin{figure}[h!]
    \dialecticaCoProductFig{}
    \caption{Product and coproduct of two objects in \(\dialLineale{L}\).}\label{fig:dialectica-objects-co-product}
  \end{figure}
  
  The terminal object is \((1 \dialarrow{\emptyset} \emptyset)\), with one element on the left, while the initial one is \((\emptyset \dialarrow{\emptyset} 1)\), with one element on the right.
\end{example}

Now, we use the monoidal operation of $L$ to define a monoidal product in $\dialLineale{L}$.
This operation can be viewed as a weaker form of ``and''~\cite{blass1992game}.
A winning strategy for Proponent in \(A \tensor B\) is a winning strategy in both \(A\) and \(B\), but, contrary to the categorical product, a winning strategy for Opponent in \(A \tensor B\) consists of a winning strategy both in \(A\) and in \(B\), for every choice of \(u \in U\) and \(v \in V\) by Proponent.
The types of \(A \tensor B\) reflect this dynamic: \(U \times V\) for the Proponent side and \(\tensortypesecond{U}{X}{V}{Y}\) for the Opponent side.

\begin{definition}[Monoidal product in \(\dialLineale{L}\)]\label{def:tensor}
  Given two objects $A = (U \dialarrow{\alpha} X)$ and $B = (V \dialarrow{\beta} Y)$ in $\dialLineale{L}$, we define their monoidal product $A \tensor B$ as the following object.
  \begin{equation*}
    A \tensor B = (U \times V  \dialarrow{\alpha  \tensor  \beta} \tensortypesecond{U}{X}{V}{Y})
  \end{equation*}
  Where \(X^{V}\) and \(Y^U\) are internal-hom objects in \(\Set\) and the function \(\alpha \tensor \beta\) is defined by the following composition.
  \begin{multline*}
    U\times V \times X^V \times Y^U \xrightarrow{\diagonal{U \times V} \times \id{X^V \times Y^U}} U \times V \times U \times V \times X^V \times Y^U \xrightarrow{\id{U \times V} \times \swap{U,V, X^{V}} \times \id{Y^{U}}} \\
    U \times V \times X^{V} \times V \times U \times Y^{U} \xrightarrow{\id{U} \times \eval_{V} \times \id{V} \times \eval_{U}} U\times X \times V \times Y\xrightarrow{\alpha \times \beta}L\times L \xrightarrow{\lintensor}L
  \end{multline*}
  where \(\diagonal{U \times V}\) is the diagonal map on \(U \times V\), \(\swap{U,V,X^{V}}\) is a permutation, and \(\eval_{U}\) and \(\eval_{V}\) are the evaluation maps in \(\Set\).
  Spelling out this definition element wise, we obtain \((\alpha \tensor \beta) (u,v,f,g) = \alpha(u,fv) \lintensor \beta(v,gu)\).

  On morphisms \((f,F) \colon A \to A'\) and \((g,G) \colon B \to B'\), we define the monoidal product as follows
  \[(f,F) \tensor (g,G) = (f \times g, F(-)g \times G(-)f) \colon A \tensor B \to A' \tensor B'\]
  where \(f \times g \colon U \times V \to U' \times V'\) and \(F(-)g \times G(-)f \colon X'^{V'}\times Y'^{U'} \to X^V \times Y^U \).
\end{definition}

\begin{proposition}\label{prop:tensor}
  The construction above induces a functor \(\tensor \colon \dialLineale{L} \times \dialLineale{L} \to \dialLineale{L}\), which is a symmetric monoidal product on \(\dialLineale{L}\).
\end{proposition}
\begin{proof}
  The object $A \tensor B = (U \times V \dialarrow{\alpha \tensor \beta}X^{V} \times Y^{U})$ is clearly an object of $\dialLineale{L}$.
  The unit is the object $\unit =(1  \dialarrow{e}1)$, which assigns to $1\times 1$ the monoidal unit $e$ of $L$.
  We need to check that the monoidal product is well defined, which means that \((f,F) \tensor (g,G) \) satisfies the condition on morphisms.
  \begin{align*}
    & \alpha \tensor \beta (u,v,(F(-)g \times G(-)f)(f',g')) \\
    & = \alpha \tensor \beta (u,v,Ff'g, Gg'f) \\
    & = \alpha(u,Ff'Gv) \lintensor \beta(v,Gg'fu) \\
    & \linorder \alpha'(fu,f'Gv) \lintensor \beta'(gv,g'fu) \\
    & = \alpha' \tensor \beta' (fu,gv,f',g') \\
    & = \alpha' \tensor \beta' ((f\times g)(u,v),f',g')
  \end{align*}
  The monoidal product is a functor as it preserves composition
  \begin{align*}
    & ((f',F') \compose (f,F)) \tensor ((g',G') \compose (g,G)) \\
    & = (f' \compose f,F \compose F') \tensor (g' \compose g,G \compose G')\\
    & = ((f' \compose f) \times (g' \compose g),FF'(-)g'g \times GG'(-)f'f)\\
    & = ((f' \times g') \compose (f \times g), (F(-)g \times G(-)f) \compose (F'(-)g' \times G'(-)f'))\\
    & = (f' \times g',F'(-)g' \times G'(-)f') \compose (f \times g,F(-)g \times G(-)f)\\
    & = ((f',F') \tensor (g',G')) \compose ((f,F) \tensor (g,G))
  \end{align*}
  and identities
  \begin{align*}
    &(\id{U},\id{X}) \tensor (\id{V},\id{Y})\\
    & = (\id{U} \times \id{V}, \id{X}(-)\id{V} \times \id{Y}(-)\id{U})\\
    & = (\id{U \times X}, \id{X^V \times Y^U})
  \end{align*}
  The associator is defined by the following isomorphisms in \(\Set\)
  \[\alpha_{A,B,C} = (\alpha_{U,V,W},A_{X,Y,Z}) \colon (A \tensor B) \tensor C \to A \tensor (B \tensor C)\]
  where \(\alpha_{U,V,W} \colon (U \times V) \times W \to U \times (V \times W)\) is the associator in \(\Set\) with the cartesian product and \(A_{X,Y,Z} \colon \tensortypesecond{U}{X}{V \times W}{(\tensortypesecond{V}{Y}{W}{Z})} \to \tensortypesecond{U \times V}{(\tensortypesecond{U}{X}{V}{Y})}{W}{Z}\) is the composition of isomorphisms in \(\Set\) given by
  \begin{multline*}
    X^{V \times W} \times (Y^W \times Z^V)^U \xrightarrow{\cong} X^{V \times W} \times (Y^{U\times W} \times Z^{U\times V})\\
    \xrightarrow{\cong} (X^{V \times W} \times Y^{U\times W}) \times Z^{U\times V} \xrightarrow{\cong} (X^{V} \times Y^{U})^W \times Z^{U\times V}
  \end{multline*}
  The unitors are defined by the following isomorphisms in \(\Set\)
  \begin{align*}
    & \lambda_A = (\lambda_U,L_X) \colon \unit \tensor A \to A && \rho_A = (\rho_U,R_X) \colon A \tensor \unit \to A
  \end{align*}
  where \(\lambda_U \colon 1 \times U \to U\) and \(\rho_U \colon U \times 1 \to U\) are the unitors in \(\Set\), and \(L_X \colon X \to 1^U \times X^1\) and \(R_X \colon X \to X^1 \times 1^U\) are the compositions of isomorphisms in \(\Set\) given by
  \begin{align*}
    &X \xrightarrow{\cong} 1 \times X \xrightarrow{\cong} 1^U \times X^1
    && X \xrightarrow{\cong} X \times 1 \xrightarrow{\cong} X^1 \times 1^U\end{align*}
  We are left to prove that the above are actually morphisms in \(\dialLineale{L}\), that they are natural isomorphisms and that they satisfy the pentagon and triangle equations~\cite{maclane}.
  The associator is a morphism because for all \(((u,v),w) \in (U\times V) \times W\) and all \((f,(g,h)) \in X^{V\times W} \times (Y^W \times Z^V)^U\)
  \begin{align*}
    & ((\alpha \tensor \beta) \tensor \gamma) (((u,v),w),A_{X,Y,Z}(f,(g,h)))\\
    & = ((\alpha \tensor \beta) \tensor \gamma) (((u,v),w),((f,g),h))\\
    & = (\alpha (u,f(v,w)) \lintensor \beta(v,g(u,w))) \lintensor \gamma(w,h(u,v))\\
    & = \alpha (u,f(v,w)) \lintensor (\beta(v,g(u,w)) \lintensor \gamma(w,h(u,v)))\\
    & = (\alpha \tensor (\beta \tensor \gamma)) ((u,(v,w)),(f,(g,h)))\\
    & = (\alpha \tensor (\beta \tensor \gamma)) (\alpha_{U,V,W}((u,v),w),(f,(g,h)))
  \end{align*}
  The unitors are morphisms because for all \(u \in U\) and all \(x \in X\)
  \begin{align*}
    & (\unit \tensor \alpha) ((\ast,u),L_{X}(x)) && (\alpha \tensor \unit) ((u,\ast),R_{X}(x)) \\
    & = (\unit \tensor \alpha) ((\ast,u),(\ast,x)) && = (\alpha \tensor \unit) ((u,\ast),(x,\ast)) \\
    & = \unit(\ast,\ast) \lintensor \alpha(u,x) && = \alpha(u,x) \lintensor \unit(\ast,\ast) \\
    & = e \lintensor \alpha(u,x) && = \alpha(u,x) \lintensor e \\
    & = \alpha(u,x) && = \alpha(u,x)
  \end{align*}
  The associator and the unitors are natural isomorphisms because they are natural isomorphisms component wise.
  The triangle and pentagon equations hold because they hold in \(\Set\).

  The symmetries \(\swap{A,B} = (\swap{U,V}, \swap{Y^{U}, X^{V}}) \colon A \tensor B \to B \tensor A\) are well-defined because their components are isomorphisms and they are natural because their components are natural in \(\Set\).
  Finally, the hexagon equations hold because they also hold in \(\Set\) with the cartesian product.
\end{proof}

\begin{definition}[Internal-hom in \(\dialLineale{L}\)]\label{def:inthom}
  Given two objects $A=(U \dialarrow{\alpha} X)$ and $B=(V \dialarrow{\beta} Y)$ in  $\dialLineale{L}$ we define their internal-hom, $\inthom{A}{B}$, as follows:
  \begin{equation*}
    \inthom{A}{B} = \inthomtypefirst{U}{X}{V}{Y} \dialarrow{\inthom{\alpha}{\beta}} U\times Y
  \end{equation*}
  The function \(\inthom{\alpha}{\beta}\) is defined by the following composition.
  \begin{multline*}
    \inthomtypefirst{U}{X}{V}{Y} \times U \times Y \xrightarrow{\id{\inthomtypefirst{U}{X}{V}{Y}} \times \diagonal{U \times Y}} \inthomtypefirst{U}{X}{V}{Y} \times U \times Y \times U \times Y \xrightarrow{\swap{}}\\
    U \times Y \times X^{Y} \times U \times V^{U} \times Y \xrightarrow{\id{U} \times \eval \times \eval \times \id{Y}} U \times X \times V \times Y \xrightarrow{\alpha \times \beta} L \times L \xrightarrow{\lininthom} L
  \end{multline*}
  Spelling out this definition element wise, we obtain $\inthom{\alpha}{\beta} (f,F,u,y)= \alpha(u,Fy) \lininthom \beta(fu, y)$.
\end{definition}

\begin{example}\label{ex:dialectica-objects-tensor-hom}
  The monoidal product and internal-hom of the objects \((U \dialarrow{\alpha} X)\) and \((V \dialarrow{\beta} Y)\) defined in \Cref{ex:dialectica-objects-morphisms} are below and in \Cref{fig:dialectica-objects-tensor-hom}.
  The sets \(X^{V} \cong \{x,x'\}\), \(Y^{U} \cong \{!\}\), \(V^{U} \cong \{!\}\) and \(X^{Y} \cong \{x,x'\}\) contain at most two elements.
  \begin{align*}
    \alpha \tensor \beta (u,v,x,!) & = a \lintensor b & \inthom{\alpha}{\beta} (!,x,u,y) & = a \lininthom b \\
    \alpha \tensor \beta (u',v,x,!) & = a' \lintensor b & \inthom{\alpha}{\beta} (!,x',u,y) & = b \\
    \alpha \tensor \beta (u,v,x',!) & = b & \inthom{\alpha}{\beta} (!,x,u',y) & = a' \lininthom b \\
    \alpha \tensor \beta (u',v,x',!) & = b & \inthom{\alpha}{\beta} (!,x',u',y) & = b 
  \end{align*}
    The monoidal unit is the object \(1 \dialarrow{e} 1\), with two elements that are not related to each other, i.e.~related with unit weight \(e\).
  \begin{figure}[h]
    \dialecticaTensorHomFig{}
    \caption{Tensor product and internal-hom of two objects in \(\dialLineale{L}\).}\label{fig:dialectica-objects-tensor-hom}
  \end{figure}

\end{example}

\begin{proposition}\label{prop:inthom}
  The construction above induces an internal-hom functor $\inthom{-}{-} \colon \dialLineale{L}\op \times \dialLineale{L} \to \dialLineale{L}$.
\end{proposition}
\begin{proof}
  The object \(\inthom{A}{B} = \inthomtypefirst{U}{X}{V}{Y} \dialarrow{\inthom{\alpha}{\beta}} U\times Y\) is clearly an object of \(\dialLineale{L}\). On morphisms \((f,F) \colon A' \to A\) and \((g,G) \colon B \to B'\) in \(\dialLineale{L}\), the internal-hom can be defined as
  \begin{equation*}
    \inthom{(f,F)}{(g,G)} = (g(-)f \times F(-)G, f \times G) \colon \inthom{A}{B} \to \inthom{A'}{B'}
  \end{equation*}
  where \(g(-)f \times F(-)G \colon \inthomtypefirst{U}{X}{V}{Y} \to \inthomtypefirst{U'}{X'}{V'}{Y'}\) and \(f \times G \colon U' \times Y' \to U \times Y\).
  We need to check that the internal-hom is well defined, which means that \(\inthom{(f,F)}{(g,G)}\) needs to satisfy the condition on morphisms.
  For all \((h,H) \in \inthomtypefirst{U}{X}{V}{Y}\) and all \((u',y') \in U' \times Y'\)
  \begin{align*}
    & \inthom{\alpha}{\beta}(h,H,(f \times G)(u',y')) \\
    & = \inthom{\alpha}{\beta}(h,H,fu', Gy') \\
    & = \alpha(fu',HGy') \lininthom \beta(hfu',Gy') \\
    & \linorder \alpha'(u',FHGy') \lininthom \beta'(ghfu',y')\\
    & = \inthom{\alpha'}{\beta'}(ghf,FHG,u',y') \\
    & = \inthom{\alpha'}{\beta'}((g(-)f \times F(-)G)(h,H),u',y')
  \end{align*}
  because \(\alpha'(u',FHGy') \linorder \alpha(fu',HGy')\) and \(\beta(hfu',Gy') \linorder \beta'(ghfu',y')\) as \((f,F)\) and \((g,G)\) are morphisms.
  The internal-hom is a functor as it preserves composition
  \begin{align*}
    & \inthom{(f',F') \compose (f,F)}{(g',G') \compose (g,G)} \\
    & = \inthom{(ff',F'F)}{(g'g,GG')} \\
    & = (g'g(-)ff' \times F'F(-)GG', ff' \times GG') \\
    & = ((g'(-)f' \times F'(-)G') \compose(g(-)f \times F(-)G),(f \times G) \compose (f' \times G')) \\
    & = (g'(-)f' \times F'(-)G',f'\times G') \compose (g(-)f \times F(-)G,f \times G) \\
    & = \inthom{(f',F')}{(g',G')} \compose \inthom{(f,F)}{(g,G)}
  \end{align*}
  and identities
  \begin{align*}
    &\inthom{\id{A}}{\id{B}}\\
    & = \inthom{(\id{U},\id{X})}{(\id{V},\id{Y})} \\ 
    & = (\id{V}(-)\id{U} \times \id{X}(-)\id{Y}, \id{U} \times \id{Y}) \\
    & = (\id{\inthomtypefirst{U}{X}{V}{Y}},\id{U \times Y}) \\
    & = \id{\inthom{A}{B}}
  \end{align*}
\end{proof}

The next result combines the structure on \(\dialLineale{L}\) defined so far.
Note how products and coproducts come from products and coproducts in \(\Set\), while the monoidal closed structure is lifted from the corresponding structure in the lineale \(L\).

\begin{theorem}\label{th:monclosed}
  The category  $\dialLineale{L}$ has products and coproducts as in \Cref{def:cartesian} and is a symmetric monoidal closed category with monoidal product as in \Cref{def:tensor} and internal-hom as in \Cref{def:inthom}.
\end{theorem}
\begin{proof}
  To prove the adjunction $-\tensor B \dashv \inthom{B}{-}$ we have to show that, for every objects \(A\) and \(C\) in \(\dialLineale{L}\), there is a bijection $\psi_{A,C} \colon \Hom{\dialLineale{L}}(A\tensor B,C)\cong \Hom{\dialLineale{L}}(A,\inthom{B}{C})$ that is natural in \(A\) and \(C\).

  Let \(\phi_{U,Z} \colon \Hom{\Set}(U \times V,Z) \to \Hom{\Set}(U,Z^V)\) be the natural isomorphism witnessing the adjunction between the cartesian product and the internal-hom in \(\Set\) and let \(\swap{U,V} \colon U \times V \to V \times U\) be the symmetry of the cartesian product in \(\Set\). Define the maps
  \begin{align*}
    & \psi_{A,C}(f,F) = (\productmap{\phi(f)}{\phi(\phi^{-1}(F_2) \compose \swap{U,Z})},\phi^{-1}(F_1) \compose \swap{V,Z}) \\
    & \psi_{A,C}^{-1}(g,G) = (\phi^{-1}(g_1), \productmap{\phi(G \compose \swap{Z,V})}{\phi(\phi^{-1}(g_2)\compose \swap{Z,U})})
  \end{align*}
  We can check that they are well defined.
  An element of $\Hom{\dialLineale{L}}(A\tensor B,C)$ is a pair $(f,\productmap{F_1}{F_2})$ with \(f \colon U \times V \to W\) and \(F = \productmap{F_1}{F_2} \colon Z \to \tensortypesecond{U}{X}{V}{Y}\)
  such that \(\forall (u,v) \in U \times V \ \forall z \in Z \ (\alpha \tensor \beta)(u,v,Fz) \linorder \gamma(f(u,v),z)\), which is equivalent to
  \(\alpha(u,(F_1(z))(v)) \lintensor \beta(v,(F_2(z))(u)) \linorder \gamma(f(u,v),z)\).
  On the other hand, an element of \(\Hom{\dialLineale{L}}(A,\inthom{B}{C})\) is a pair \((\productmap{g_1}{g_2},G)\) with \(g = \productmap{g_1}{g_2} \colon U \to \inthomtypefirst{V}{Y}{W}{Z}\) and \(G \colon V \times Z \to X\)
  such that \(\forall u \in U \ \forall (v,z) \in V \times Z \ \alpha(u,G(v,z)) \linorder \inthom{\beta}{\gamma}(g(u),v,z)\), which is equivalent to
  \(\alpha(u,G(v,z)) \linorder \beta(v,(g_2(u))(z)) \lininthom \gamma((g_1(u))(v),z)\).

  They are morphisms because the inequality condition for morphisms in \(\dialLineale{L}\) holds with equality.
  We check that they are inverses to each other.
  \begin{align*}
    & \psi_{A,C} \compose \psi_{A,C}^{-1}(g,G) \\
    & = (\productmap{\phi(\phi^{-1}(g_1))}{\phi(\phi^{-1}(\phi(\phi^{-1}(g_2)\compose \swap{Z,U}))\compose \swap{U,Z})}, \phi^{-1}(\phi(G \compose \swap{Z,V}))\compose\swap{V,Z}) \\
    & = (\productmap{g_1}{g_2},G) \\
    & \\
    & \psi_{A,C}^{-1}\compose \psi_{A,C}(f,F) \\
    & = (\phi^{-1}(\phi(f)),\productmap{\phi(\phi^{-1}(F_1)\compose \swap{V,Z}\compose \swap{Z,V})}{\phi(\phi^{-1}(\phi(\phi^{-1}(F_2)\compose \swap{U,Z}))\compose \swap{Z,U})}) \\
    & = (f,\productmap{F_1}{F_2})
  \end{align*}
  We check that they are natural.
  Naturality comes from naturality of \(\phi\) in \(\Set\).
\end{proof}

\section{Dialectica Petri nets}\label{sec:Petri-net-cat}

A Petri net is given by a set of places \(U\), a set of transitions \(X\), and has two relations between these two sets that specify the precondition relation \(\pre{\alpha}\) and the postcondition relation \(\post{\alpha}\).
In our case these relations will be valued in a generic lineale \(L\) and the pre- and post- conditions will be objects in \(\dialLineale{L}\).
The category of Petri nets that we consider has Petri nets as objects and is obtained by putting together two copies of \(\dialLineale{L}\) by taking a pullback in \(\Cat\): the first copy represents preconditions \(\pre{\alpha} \colon U \times X \to L\) and the other one represents postconditions \(\post{\alpha} \colon U \times X \to L\).

\begin{definition}[Category \(\dialNet{L}\)]\label{def:dialectica-PN}
Given a lineale \((L,\linorder,\lintensor,e,\lininthom)\), the category \(\dialNet{L}\) is defined by the following data.
\begin{itemize}
    \item An object is a pair \(A = (\pre{A},\post{A})\) of objects \(U \dialarrow{\pre{\alpha}} X\) and \(U \dialarrow{\post{\alpha}} X\) in \(\dialLineale{L}\).
    \item A morphism \((f,F) \colon (\pre{A},\post{A}) \to (\pre{B},\post{B})\) is both a morphism \((f,F) \colon \pre{A} \to \pre{B}\) and a\\  morphism \((f,F) \colon \post{A} \to \post{B}\) in \(\dialLineale{L}\).
\end{itemize}
\end{definition}

\begin{remark}\label{rem:nets-pullback}
  \Cref{def:dialectica-PN} can be restated more abstractly.
  The category \(\dialNet{L}\) is the pullback in \(\Cat\) of the forgetful functor \(\mathcal{U}_{L} \colon \dialLineale{L} \to \Set \times \Set\op\) along itself.
  The functor \(\mathcal{U}_{L}\) assigns to an object \((U,X,\alpha)\) the pair of sets \((U,X)\), and to a morphism \((f,F)\) the corresponding pair of functions \((f,F)\).
\end{remark}

\begin{example}\label{ex:nets-morphisms}
  The relations in \Cref{ex:dialectica-objects-morphisms} may represent the precondition relations, \((U \dialarrow{\pre{\alpha}} X)\) and \((V \dialarrow{\pre{\beta}} Y)\), of two nets \(A\) and \(B\).
  We define below two other relations, \((U \dialarrow{\post{\alpha}} X)\) and \((V \dialarrow{\post{\beta}} Y)\), that specify the postconditions of the nets \(A\) and \(B\). \Cref{fig:nets-objects-morphisms} shows the net \(A\) on the top left and bottom right, and the net \(B\) on the bottom left and top right.
  \begin{align*}
    \pre{\alpha}(u,x) & = a & \pre\alpha(u',x) & = a'& \post{\alpha}(u,x) & = e & \post{\alpha}(u',x) & = e \\
    \pre\alpha(u,x') & = e & \pre\alpha(u',x') & = e & \post{\alpha}(u,x') & = c & \post{\alpha}(u',x') & = e \\
    \pre\beta(v,y) & = b & & & \post{\beta}(v,y) & = e &&
  \end{align*}
  The morphism \((f,F)\) defined in \Cref{ex:dialectica-objects-morphisms} is also a morphism \((f,F) \colon A \to B\) of nets because it is a morphism \((f,F) \colon (U \dialarrow{\post{\alpha}} X) \to (V \dialarrow{\post{\beta}} Y)\), see \Cref{fig:nets-objects-morphisms}, left.
  The morphism \((g,G)\) in the same example is a morphism of nets whenever \(e \linorder c\), see \Cref{fig:nets-objects-morphisms}, right.
  This condition, for the lineales like \(2\), \(3\) or \(\Naturals\) where \(e\) is the smallest element, is satisfied only for \(c = e\) so the elements \(u\) and \(x'\) cannot be related.
  \begin{figure}[h!]
    \netsMorphismFig{}
    \caption{Two net morphisms, if \(a \linorder b\) (left), and if \(b \linorder a,a'\) and \(e \linorder c\) (right).\label{fig:nets-objects-morphisms}}
  \end{figure}
\end{example}

% \begin{example}
% Let us consider some examples of morphisms in \(\dialNet{\Naturals}\).
% \begin{itemize}
%   \item Simulations regarding $\alpha$ and $\beta$: a morphism \((f,F) \colon A \to B\) in this category may represent the fact that the Petri net \(A\) is a simulation of the net \(B\) as the conditions on the morphisms in \(\dialLineale{\Naturals}\) ensure that the preconditions and the postconditions of \(B\) are `smaller or equal' than those of \(A\): \(\forall u \in U \ \forall y \in Y \pre{\alpha}(u,Fy) \geq \pre{\beta}(fu,y) \land \post{\alpha}(u,Fy) \geq \post{\beta}(fu,y)\).
%         \Cref{fig:morphisms} (left) illustrates this type of morphisms.
%   \item Refinements regarding $f$: in \Cref{fig:morphisms} (right), the Petri net \(B\)  represents a refinement of the net \(A\) obtained by increasing its number of components (by adding a place).
% \end{itemize}
% \begin{figure}[ht]
% \begin{subfigure}[ht]{0.45\linewidth}
% \centering
% \includegraphics[width=.8\linewidth]{simulation.pdf}
% \end{subfigure}
% \hfill
% \begin{subfigure}[ht]{0.45\linewidth}
% \centering
% \includegraphics[width=.8\linewidth]{refinement.pdf}
% \end{subfigure}
% \caption{Examples of morphisms in \(\dialNet{\Naturals}\): simulations (left) and refinements (right).}\label{fig:morphisms}
% \end{figure}
% \end{example}

The structure of \(\dialLineale{L}\) defines analogous structure in \(\dialNet{L}\).
\begin{definition}[Structure of \(\dialNet{L}\)]
The category \(\dialNet{L}\) inherits the structure of \(\dialLineale{L}\). All the connectives are defined componentwise:
\begin{itemize}
    \item \(A \tensor B = (\pre{A} \tensor \pre{B}, \post{A} \tensor \post{B})\).
    \item \(\inthom{A}{B} = (\inthom{\pre{A}}{\pre{B}}, \inthom{\post{A}}{\post{B}})\).
    \item \(A \product B = (\pre{A} \product \pre{B}, \post{A} \product \post{B})\).
    \item \(A \coproduct B = (\pre{A} \coproduct \pre{B}, \post{A} \coproduct \post{B})\).
\end{itemize}
\end{definition}
Examples of Petri nets modelled in this category are in the next section, where we will show how, with the possibility of changing the lineale, we can encompass different kinds of nets.

\begin{example}\label{ex:nets-co-product}
  \Cref{fig:nets-co-product} shows the product and coproduct of the nets defined in \Cref{ex:nets-morphisms}.
  These are obtained by computing their product and coproduct component-wise, as shown in \Cref{ex:dialectica-objects-co-product}.
  \Cref{fig:nets-terminal-initial-unit} shows the terminal and initial nets, which are the units for the product and coproduct respectively.
  \begin{figure}[h!]
    \netsCoProductFig{}
    \caption{Product (left) and coproduct (right) of the nets \(A\) and \(B\).\label{fig:nets-co-product}}
  \end{figure}
\end{example}

\begin{figure}[h!]
  \netsTerminalInitialUnit{}
  \caption{Terminal net (left), initial net (centre) and unit net (right).\label{fig:nets-terminal-initial-unit}}
\end{figure}

\begin{example}\label{ex:nets-tensor-hom}
  \Cref{fig:nets-tensor-hom} shows the monoidal product and internal-hom of the nets defined in \Cref{ex:nets-morphisms}.
  These are obtained by computing their monoidal product and internal-hom component-wise, as shown in \Cref{ex:dialectica-objects-tensor-hom}.
  \Cref{fig:nets-terminal-initial-unit} shows the monoidal unit net.
  \begin{figure}[h!]
    \netsTensorHomFig{}
    \caption{Monoidal product (left) and internal-hom (right) of the nets \(A\) and \(B\).\label{fig:nets-tensor-hom}}
  \end{figure}
\end{example}

\section{Changing the logic of Petri nets}\label{sec:changing-logic}

This section studies in detail the examples mentioned in \Cref{sec:lineales-examples} to showcase the different interpretations of the arcs that we can achieve by just changing the lineale \(L\) in the construction of the category \(\dialNet{L}\).
\Cref{sec:functors-from-morphisms} concludes by making the constructions in \Cref{sec:intermidiate-cat} and \Cref{sec:Petri-net-cat} functorial: for a morphism of lineales \(h \colon L \to L'\), specifying how to change the logic, we construct a functor \(\dialLineale{L} \to \dialLineale{L'}\) that preserves the structure of these categories.

\subsection{Elementary Petri nets: \(L=2\)}
By considering relations with values on this lineale, which correspond to ordinary relations, we obtain the elementary Petri nets, those where pre- and post-conditions only say whether or not a place is a pre- or post-condition for a transition~\cite{rozenberg1996elementary}.

\subsection{Kleene Petri nets: \(L=3\)}
Thanks to \Cref{lem:kleene-lineale}, we can define the Dialectica construction over \((3,\leq,\lintensor,1,\lininthom)\) and Petri nets with weights in \(3\) accordingly.

We take as motivating example the model of the chemical reactions regulating the circadian clock of \textit{Synechococcus Elongatus}~\cite{synechococcus} that is composed of two successive phosphorylations and two successive dephosphorilations (which are the transitions labelled with ``p'' and ``d'', respectively, in \Cref{fig:synechococcus}).
There is experimental evidence~\cite{axmann2007minimal} for the existence of further feedback loops in this model.
However, the precise underlying mechanism is still unknown.
We can take into account these unknowns in our model for Petri nets by adding arcs with \(u\) weight (presence and absence are represented by \(1\) and \(0\) respectively).
The Petri net in \Cref{fig:synechococcus} shows the values of the pre- and post- conditions relations as weights on the arcs.
Note that the arcs not shown are those with value \(0\).

\begin{figure}[h!]
    \centering
    \begin{tikzpicture}[xscale=1.5]
        \node[place,tokens=0,label=center:P](p1)at(0,4){};
        \node[place,tokens=0,label=below:KaiA](ka1)at(0,2){};
        \node[place,tokens=0,label=above:KaiA](ka2)at(2,6){};
        \node[place,tokens=0,label=above:KaiBC+P](kbcp)at(2,4){};
        \node[place,tokens=0,label=below left:KaiABC+P](kabcp)at(2,2){};
        \node[place,tokens=0,label=above:KaiB](kb1)at(4,6){};
        \node[place,tokens=0,label=center:P](p2)at(4,5){};
        \node[place,tokens=0,label=below:KaiAC](kac)at(4,4){};
        \node[place,tokens=0,label=above:KaiAC+P](kacp)at(4,2){};
        \node[place,tokens=0,label=right:KaiB](kb2)at(4,1){};
        \node[place,tokens=0,label=center:P](p4)at(4,0){};
        \node[place,tokens=0,label=center:P](p3)at(6,3){};
        \node[transition,label=center:d](dp1)at(1,3){}
        edge[post]node[above]{$1$}(p1)
        edge[post]node[above]{$1$}(kbcp)
        edge[post]node[above]{$1$}(ka1)
        edge[pre]node[above]{$1$}(kabcp)
        edge[pre]node[below]{$u$}(kac);
        \node[transition,label=center:d](dp2)at(3,5){}
        edge[post]node[above]{$1$}(kb1)
        edge[post]node[above]{$1$}(p2)
        edge[post]node[above]{$1$}(kac)
        edge[pre]node[above]{$1$}(kbcp)
        edge[pre]node[above]{$1$}(ka2);
        \node[transition,label=center:p](ph1)at(5,3){}
        edge[post]node[above]{$1$}(kacp)
        edge[pre]node[above]{$1$}(p3)
        edge[pre]node[above]{$1$}(kac);
        \node[transition,label=center:p](ph2)at(3,1){}
        edge[post]node[above]{$1$}(kabcp)
        edge[pre]node[above]{$1$}(kacp)
        edge[pre]node[above]{$1$}(kb2)
        edge[pre]node[above]{$1$}(p4)
        edge[pre]node[right]{$u$}(kbcp);
    \end{tikzpicture}
      \captionsetup{width=.9\linewidth}
    \caption{Petri net representing the chemical reaction network regulating the circadian clock of Synechococcus Elongatus. Present and unknown relations are labelled by \(1\) and \(u\), respectively.}\label{fig:synechococcus}
\end{figure}

\subsection{Multirelation Petri nets: \(L=\Naturals\)}

As every chemical reaction, the one to obtain water from oxygen and hydrogen needs stoichiometric coefficients to be represented properly. Multirelations take these into account, as shown in \Cref{fig:PNwater}.
\begin{figure}[h!]
    \centering
    \begin{tikzpicture}[xscale=1.5]
        \node[place,tokens=0,label=center:H$_2$](h2)at(0,2){};
        \node[place,tokens=0,label=center:O$_2$](o2)at(0,0){};
        \node[place,tokens=0,label=center:H$_2$O](h2o)at(2,1){};
        \node[transition](t)at(1,1){}
        edge[post]node[above]{$2$}(h2o)
        edge[pre]node[above]{$2$}(h2)
        edge[pre]node[above]{$1$}(o2);
    \end{tikzpicture}
    \caption{Petri net representing the chemical reaction 2H$_2$ + O$_2$ $\to$ 2H$_2$O.}\label{fig:PNwater}
\end{figure}

\subsection{Integers Petri nets: \(L=\Integers\)}

  Empirical systems often need to locally reverse the logic of preconditions to express that the presence of tokens in a given place ``disables'' a transition.
  Several different concepts of inhibitor arcs can be modelled by Petri nets including the ``threshold inhibitor arc''.
  Reaction inhibitors in chemistry illustrate the situation: in \Cref{fig:inihibitorarc} chemical reaction $r$ will not take place if the amount of substance I exceeds 3, a condition that is expressed by its inverse $-3$.

\begin{figure}[h!]
    \centering
    \begin{tikzpicture}[xscale=1.5]
        \node[place,tokens=0,label=center:S$_1$](s1)at(0,2){};
        \node[place,tokens=0,label=center:S$_2$](s2)at(0,0){};
        \node[place,tokens=0,label=center:S$_3$](s3)at(2,1){};
        \node[place,tokens=0,label=center:I](i)at(1,2.5){};
        \node[transition](t)at(1,1){$r$}
        edge[pre]node[above]{$2$}(s1)
        edge[pre]node[right]{$-3$}(i)
        edge[pre]node[above]{$2$}(s2)
        edge[post]node[below]{$1$}(s3);
    \end{tikzpicture}
      \captionsetup{width=.9\linewidth}
    \caption{Petri net representation of the chemical reaction S$_1$ + S$_2$ $\to$ S$_3$. The inhibitor arc is labelled by -3, expressing that 3 is the minimum amount of substance I that prevents $r$ from taking place.}\label{fig:inihibitorarc}
\end{figure}

\subsection{Probabilistic Petri nets: \(L=[0,1]\)}

The SIR (Susceptible, Infectious, Recovered) model is a simple compartmental model for infectious diseases. A susceptible individual has a contact with an infectious individual with probability \(p_c\) and, after the contact, it can be infected with probability \(p_I\), or remain susceptible with probability \(1-p_I\). On the other hand, an infectious individual can recover with probability \(p_R\) or remain infectious with probability \(1-p_R\).
This setting can be represented with a Petri net where the relations between places and transitions are valued in \([0,1]\) (\Cref{fig:PNsir}).
\begin{figure}[h!]
    \centering
    \begin{tikzpicture}[xscale=1.5]
        \node[place,tokens=0,label=center:S](S)at(0,0){};
        \node[place,tokens=0,label=center:I](I)at(2,0){};
        \node[place,tokens=0,label=center:R](R)at(4,1){};
        \node[transition,label=center:c](cont)at(1,0){}
        edge[post,out=45,in=135]node[above]{$p_I$}(I)
        edge[post,out=-135,in=-45]node[below]{$1-p_I$}(S)
        edge[pre,out=135,in=45]node[above]{$p_c$}(S)
        edge[pre,out=-45,in=-135]node[below]{$1$}(I);
        \node[transition,label=center:r](rec)at(3,1){}
        edge[post]node[above]{$1$}(R)
        edge[pre]node[above]{$p_R$}(I);
        \node[transition,label=center:i](inf)at(3,-1){}
        edge[post,out=90,in=0]node[above]{$1$}(I)
        edge[pre,out=180,in=-90]node[below]{$1-p_R$}(I);
    \end{tikzpicture}
    \caption{Petri net representing the SIR model.}\label{fig:PNsir}
\end{figure}

\subsection{Product of lineales}
  There is a dual situation to inhibition in chemistry, namely, catalysis. A catalyst is a substance that increases the reaction rate without being consumed by the reaction.
  The presence of a substance S in a chemical reaction might then play one of three roles: reactant/product, inhibitor, or catalyst. We claim that the product of the lineales $\Reals^+$ and $\Integers$ has enough expressive power to model reaction rates in the presence of both inhibitors and catalysts. In \Cref{fig:product} pairs of the form $(r,0)$, state that those substances are not inhibitors nor catalysts, and $r$ is the rate at which a substance is consumed or produced.
  The negative number in the label $(r_4,-3)$ expresses that I is an inhibitor of reaction $r$, and -3 the minimum amount of I required to slow down the reaction by the rate $r_4$.  Finally, the label $(r_5,5)$ indicates that C is a catalyst and 5 is the minimum amount of C required to increase the reaction rate by $r_5$.

\begin{figure}[h!]
    \centering
    \begin{tikzpicture}[xscale=1.5]
        \node[place,tokens=0,label=center:S$_1$](s1)at(-1,3){};
        \node[place,tokens=0,label=center:S$_2$](s2)at(-1,1){};
        \node[place,tokens=0,label=center:S$_3$](s3)at(3,2){};
        \node[place,tokens=0,label=center:I](i)at(1,3.5){};
        \node[place,tokens=0,label=center:C](c)at(1,0){};
        \node[transition](t)at(1,2){$r$}
        edge[pre]node[above]{$(r_1,0)$}(s1)
        edge[pre]node[right]{$(r_4,-3)$}(i)
        edge[pre]node[above]{$(r_2,0)$}(s2)
        edge[pre]node[left]{$(r_5,5)$}(c)
        edge[post]node[below]{$(r_3,0)$}(s3);
    \end{tikzpicture}
      \captionsetup{width=.9\linewidth}
    \caption{Petri net representation of reaction rates for the chemical reaction S$_1$ + S$_2$ $\to$ S$_3$ in the presence of an inhibitor I and a catalyst C. Labels are pairs $(r,z)$ where $z$ states the role of the substance as reactant/product (zero), inhibitor (negative integers), and catalysts (positive integers); and $r$ the rate at which the substance is consumed/produced (if $z=0$), or at which the reaction rate increases ($z>0$) or is slowed down (if $z<0$).}\label{fig:product}
\end{figure}

\subsection{Changing the logic, functorially}\label{sec:functors-from-morphisms}

This section concludes by extending the constructions \(\dialLineale{L}\) and \(\dialNet{L}\) to functors.
Each morphism of lineales \(h \colon L \to L'\) determines a functor \(\dialLineale{L} \to \dialLineale{L'}\) that maps \(L\)-valued relations to \(L'\)-valued relations and, thus, changes their logic.
Similarly, we obtain a functor \(\dialNet{L} \to \dialNet{L'}\) that maps \(L\)-valued nets to \(L'\)-valued nets.

Consider the category \(\LinCat\) of symmetric monoidal closed categories with products and coproducts, where morphisms are lax monoidal functors preserving the products and coproducts.
\Cref{th:monclosed} and \Cref{prop:dialectica-co-products} show that, for every lineale \(L\), \(\dialLineale{L}\) is an object of \(\LinCat\).
The next proposition shows that this construction extends to a functor preserving the closed monoidal structure, products and coproducts.

\begin{proposition}\label{prop:dial-lineale-functor}
  There is a functor \(\dialLineale{(-)} \colon \Lin \to \LinCat\) defined on objects \(L\) by the construction \(\dialLineale{L}\) in \Cref{sec:intermidiate-cat}.
\end{proposition}
\begin{proof}
  For a morphism of lineales \(h \colon (L, \linorder, \lintensor, e, \lininthom) \to (L', \linorder', \lintensor', e', \lininthom')\), define the action of \(\dialLineale{(-)}\) on it as a functor \(\dialLineale{h} \colon \dialLineale{L} \to \dialLineale{L'}\).
  On an object \(A = U \dialarrow{\alpha} X\), it acts by postcomposing with \(h\): \(\dialLineale{h}( U \dialarrow{\alpha} X) =  U \dialarrow{h \compose \alpha} X\).
  On morphisms, it acts as the identity: \(\dialLineale{h}(f,F) = (f,F)\).
  This is well-defined because \(h\) is monotone:
  \[h \compose \alpha \compose (\id{U} \times F) \linorder h \compose \beta \compose (f \times \id{Y}) \ .\]
  \(\dialLineale{h}\) is a functor as it trivially preserves compositions and identities.
  We check that it preserves products and coproducts.
  \begin{align*}
    & \dialLineale{h}(A \product B) && \dialLineale{h}(A \coproduct B)\\
    & = U \times V \dialarrow{h \compose (\alpha \product \beta)} X + Y && = U + V \dialarrow{h \compose (\alpha \coproduct \beta)} X \times Y\\
    & = U \times V \dialarrow{h \compose \coproductmap{\alpha \compose \discard_{V}}{\beta \compose \discard_{U}}} X + Y && = U + V \dialarrow{h \compose \coproductmap{\alpha \compose \discard_{Y}}{\beta \compose \discard_{X}}} X \times Y\\
    & = U \times V \dialarrow{\coproductmap{h \compose \alpha \compose \discard_{V}}{h \compose \beta \compose \discard_{U}}} X + Y && = U + V \dialarrow{\coproductmap{h \compose \alpha \compose \discard_{Y}}{h \compose \beta \compose \discard_{X}}} X \times Y\\
    & = U \times V \dialarrow{(h \compose \alpha) \product (h \compose \beta)} X + Y && = U + V \dialarrow{(h \compose \alpha) \coproduct (h \compose \beta)} X \times Y\\
    & = \dialLineale{h}(A) \product \dialLineale{h}(B) && = \dialLineale{h}(A) \coproduct \dialLineale{h}(B)
  \end{align*}
  Finally, we show that \(\dialLineale{h}\) is a lax monoidal functor.
  Since \(h\) is a lax monoid homomorphism, the pair \(\mu_{A,B} = (\id{U \times V}, \id{X^{V} \times Y^{U}}) \colon \dialLineale{h}(A) \tensor \dialLineale{h}(B) \to \dialLineale{h}(A \tensor B)\) is a morphism in \(\dialLineale{L'}\).
  \[\dialLineale{h}(A) \tensor \dialLineale{h}(B) = U \times V \dialarrow{(h \compose \alpha) \tensor (h \compose \beta)} X^{V} \times Y^{U} \linorder  U \times V \dialarrow{h \compose (\alpha \tensor \beta)} X^{V} \times Y^{U} = \dialLineale{h}(A \tensor B)\]
  Similarly, \(\varepsilon = (\id{1},\id{1}) \colon I \to \dialLineale{h}(I)\) is a morphism in \(\dialLineale{L'}\).
  \[I = 1 \dialarrow{e'} 1 \linorder 1 \dialarrow{h \compose e} = \dialLineale{h}(I)\]
  Naturality, associativity and unitality of \(\mu_{A,B}\) are easy to check as its components are identities and \(\dialLineale{h}\) is the identity on morphisms.
  Finally, \(\dialLineale{h}\) is symmetric because it is the identity on morphisms.
\end{proof}

Every category constructed by the functor in \Cref{prop:dial-lineale-functor} determines a functor \(\mathcal{U}_{L} \colon \dialLineale{L} \to \Set \times \Set\op\) that forgets the relations \(U \dialarrow{\alpha} X\) and only keeps the underlying sets \((U,X)\).
Consider the image of the functor \(\dialLineale{(-)}\), viewed as a subcategory \(\DialCat \hookrightarrow \LinCat\).

\begin{lemma}\label{lem:forgetful-functors}
  There is a functor \(\mathcal{U} \colon \DialCat \to \DialCat/(\Set \times \Set\op)\) that assigns the forgetful functor \(\mathcal{U}_{L} \colon \dialLineale{L} \to \Set \times \Set\op\) to each category \(\dialLineale{L}\).
\end{lemma}
\begin{proof}
  For a functor \(H \colon \dialLineale{L} \to \dialLineale{L'}\), we define \(\mathcal{U}(H) = H\).
  This assignment is functorial as every such functor \(H\) commutes with the forgetful functors: on objects \(A = U \dialarrow{\alpha} X\), \(\mathcal{U}_{L'}(H(A)) = (U,X) = \mathcal{U}_{L}(A)\); on morphisms \((f,F) \colon U \dialarrow{\alpha} X \to V \dialarrow{\beta} Y\), \(\mathcal{U}_{L'}(H(f,F)) = (f,F) = \mathcal{U}_{L}(f,F)\).
\end{proof}

As mentioned in \Cref{rem:nets-pullback}, the category \(\dialNet{L}\) is a pullback of the functor \(\mathcal{U}_{L} \colon \dialLineale{L} \to \Set \times \Set\op\) along itself.
This construction is functorial.

\begin{proposition}\label{prop:ker-functor}
  There is a functor \(\mathsf{Ker} \colon \Cat/(\Set \times \Set\op) \to \Cat\) that, for a functor \(F \colon \namedcat{C} \to \Set \times \Set\op\), constructs the pullback of \(F\) along itself.
\end{proposition}
\begin{proof}
  For a functor \(F \colon \namedcat{C} \to \Set \times \Set\op\), construct the pullback \(\hat{\namedcat{C}}\) of \(F\) along itself (in \(\Cat\)) and define \(\mathsf{Ker}(F) = \hat{\namedcat{C}}\).
  Let \(F \colon \namedcat{C} \to \Set \times \Set\op\) and \(G \colon \namedcat{D} \to \Set \times \Set\op\) be objects in \(\Cat/(\Set \times \Set\op)\).
  For a functor \(H \colon \namedcat{C} \to \namedcat{D}\) such that \(G \compose H = F\), define \(\mathsf{Ker}(H)\) as the unique functor \(\hat{\namedcat{C}} \to \hat{\namedcat{D}}\) given by the universal property of the pullback defining \(\hat{\namedcat{D}}\).
  This defines a functor as it is the composition of the diagonal and product functors in \(\Cat/(\Set \times \Set\op)\) with the forgetful functor \(\Cat/(\Set \times \Set\op) \to \Cat\).
\end{proof}

\begin{corollary}\label{cor:dialnet-functorial}
  There is a functor \(\dialNet{(-)} \colon \Lin \to \Cat\) defined on objects \(L\) by the construction \(\dialLineale{L}\) in \Cref{sec:Petri-net-cat}.
\end{corollary}
\begin{proof}
  The functor \(\dialNet{(-)}\) is the composition below of functors from \Cref{prop:dial-lineale-functor,prop:ker-functor,lem:forgetful-functors}.
  \[\dialNet{(-)} \colon \Lin \xrightarrow{\dialLineale{(-)}} \DialCat \xrightarrow{\mathcal{U}} \DialCat/(\Set \times \Set\op) \hookrightarrow \Cat/(\Set \times \Set\op) \xrightarrow{\mathsf{Ker}} \Cat\]
\end{proof}

\begin{example}
  Consider the morphism of lineales \(h \colon \Naturals \to 2\) in \Cref{ex:lineale-morph}.
  By \Cref{cor:dialnet-functorial}, we may construct a functor \(\dialNet{h} \colon \dialNet{\Naturals} \to \dialNet{2}\) that forgets the amount of tokens consumed or produced by transitions and only keeps which places appear as pre and post conditions for transitions.
  For example, the Petri net in \Cref{fig:PNwater} would be mapped to a Petri net without labels, i.e.\ where all the arcs are labelled by \(1\).
\end{example}

\section{Conclusions and further work}

% Most of the recent work on Petri nets focuses on the unfoldings of a single net.
Our work follows the lines of recent work on categorical compositionality of Petri nets: they can be composed along shared places~\cite{baez2020open} or along shared transitions~\cite{pawel2014reachability}.
We have explored composing nets with some of the connectives typical of linear logic: products, coproducts, tensor and internal hom.
We have presented a categorical model for Petri nets that focuses on the diverse nature of network relations.
This approach allows the use of Petri nets with different kinds of transitions (different kinds of labels in their graphs), while maintaining their composionality.

Our model can handle different kinds of transition whose labels can be represented as a lineale (a poset version of a symmetric monoidal closed category).  Several sets of labels, from those often used in empirical data modelling, can be endowed with the structure of a lineale, including: stoichiometric coefficients in chemical reaction networks ($L=\Naturals$), reaction rates ($L=\Reals^+$), inhibitor arcs ($L=\Integers$), gene interactions ($L=\{0,1\}$), unknown or incomplete data ($L=\{-1,0,1\}$), and probabilities ($L= [0,1]$).

The structure of the lineale is lifted to the category $\dialLineale{L}$ from which the category $\dialNet{L}$ of Petri nets is built.
Both $\dialLineale{L}$ and $\dialNet{L}$ are symmetric monoidal closed categories with finite products and coproducts, providing a compositional way to put together smaller nets into bigger ones, making sure that morphisms between the component nets also assemble into morphisms between the resulting nets.
Both these constructions are functorial.
% making sure that algebraic properties of the components are preserved in the resulting net.

The category $\dialNet{L}$ is a model for weighted and directed bipartite relations and therefore we anticipate applications of the compositionality of $\dialNet{L}$ in the broader context of directed bipartite graphs~\cite{Imrich2000,HAMMACK20091114}, in particular, for their applications to real-world networks.
For instance, the labelled wiring of these graphs is key to the empirical analysis of metabolic networks, where the metabolism of an organism is studied in terms of the concurrence of smaller functional subnets called modules. We wonder whether our formal connectives may assist in the reconstruction and understanding of the whole metabolism in terms of the concurrence of the modules.

There is much more work to be done still. Both in the applications we are pursuing and on the theory of Dialectica Petri nets.
On the theory side notions of behaviour (token game) should be investigated and on the practical side we have still to investigate how the implemented systems can be modified to deal with our nets.
Moreover, our exposition fixed the base category to be \(\Set\) and only varied the lineale parameter \(L\).
The possibility of changing the base category remains to be investigated.
In particular, in~\cite{dePaiva1991multirelations}, the author shows how to construct a model of Intuitionistic Linear Logic with a modality ``!" on the category \(\dialLinealeCat{L}{C}\), provided the base category \(\namedcat{C}\) is cartesian closed and has free monoids (such as \(\Set\)). The reader interested in this construction and its interpretation in terms of Petri nets is referred to~\cite{dePaiva1991multirelations}. In \(\dialLinealeCat{\Naturals}{C}\), the construction works for a fixed lineale \((\Naturals,\geq,+,0,\lininthom)\) and is parametrized by the category \(\namedcat{C}\). In light of our present work, this opens several avenues for further research. For instance, one could explore building a modality ``!" in \(\dialLineale{L}\) for any lineale \(L\), or generalizing the \(\dialLinealeCat{\Naturals}{C}\) construction to lineales other than \((\Naturals,\geq,+,0,\lininthom)\) using our approach.  

Finally, we leave for future work the investigation of the categorical properties of the constructions \(\dialLineale{(-)}\) and \(\dialNet{(-)}\), including their potential connections with fibrations and their double-categorical aspects.
We also leave for future work the investigation of categorical properties of the constructions \(\dialLineale{(-)}\) and \(\dialNet{(-)}\), their connections with fibrations or their double categorical aspects.

Dialectica Petri nets share some of the pros and cons of other Linear Logic based nets.
As far as we know no one has investigated Differential Linear Logic Petri nets, yet (see~\cite{EHRHARD} for relating differential interaction nets to the $\pi$-calculus). We wonder if this would make the exchange of information with modellers somewhat easier. Finally we would like to investigate whether we could code our nets using Catlab \url{https://github.com/AlgebraicJulia/Catlab.jl}, a framework for  computational category theory, written in the Julia language and already used for other styles of Petri nets.

\paragraph*{Acknowledgements.}
We would like to  thank the ACT Adjoint School 2020 and its organisers for first setting up our collaboration.
We also thank Jade Master, Xiaoyan Li and Eigil Rischel for their initial work in this project.

\bibliography{references}

\newcommand{\etalchar}[1]{$^{#1}$}
\begin{thebibliography}{SNMM11}

\bibitem[ALH07]{axmann2007minimal}
Ilka~M Axmann, Stefan Legewie, and Hanspeter Herzel.
\newblock A minimal circadian clock model.
\newblock In {\em Genome Informatics 2007: Genome Informatics Series Vol. 18},
  pages 54--64. World Scientific, 2007.

\bibitem[BG90]{brown1990}
Carolyn Brown and Doug Gurr.
\newblock A categorical linear framework for {P}etri nets.
\newblock In {\em [1990] Proceedings. Fifth Annual IEEE Symposium on Logic in
  Computer Science}, pages 208--218. IEEE, 1990.

\bibitem[BG95]{brown1995}
Carolyn Brown and Doug Gurr.
\newblock A categorical linear framework for {P}etri nets.
\newblock {\em Information and Computation}, 122(2):268 -- 285, 1995.

\bibitem[BGdP91]{brown1991}
Carolyn Brown, Doug Gurr, and Valeria de~Paiva.
\newblock A {L}inear {S}pecification {L}anguage for {P}etri {N}ets ({A}arhus
  {T}echnical {R}eport), 1991.

\bibitem[BGMS21]{baez2021categories}
John~C Baez, Fabrizio Genovese, Jade Master, and Michael Shulman.
\newblock Categories of nets.
\newblock {\em arXiv preprint arXiv:2101.04238}, 2021.

\bibitem[BJ77]{belnap1977useful}
Nuel~D Belnap~Jr.
\newblock A useful four-valued logic.
\newblock In {\em Modern uses of multiple-valued logic}, pages 5--37. Springer,
  1977.

\bibitem[Bla92]{blass1992game}
Andreas Blass.
\newblock A game semantics for linear logic.
\newblock {\em Annals of Pure and Applied logic}, 56(1-3):183--220, 1992.

\bibitem[BM20]{baez2020open}
John~C Baez and Jade Master.
\newblock Open {P}etri nets.
\newblock {\em Mathematical Structures in Computer Science}, 30(3):314--341,
  2020.

\bibitem[BMMS01]{bruni2001functorial}
Roberto Bruni, Jos{\'e} Meseguer, Ugo Montanari, and Vladimiro Sassone.
\newblock Functorial models for {P}etri nets.
\newblock {\em Information and Computation}, 170(2):207--236, 2001.

\bibitem[Bor94]{Borceux_1994}
Francis Borceux.
\newblock {\em Handbook of Categorical Algebra}.
\newblock Encyclopedia of Mathematics and its Applications. Cambridge
  University Press, 1994.

\bibitem[dP89a]{dePaiva1989:AMS}
Valeria de~Paiva.
\newblock The {D}ialectica {C}ategories.
\newblock In {\em Categories in Computer Science and Logic: Proceedings of the
  AMS-IMS-SIAM Joint Summer Research Conference Held June 14-20, 1987 with
  Support from the National Science Foundation}, volume~92, page~47. American
  Mathematical Society, 1989.

\bibitem[dP89b]{dePaiva1989:CTCS}
Valeria de~Paiva.
\newblock A dialectica-like model of linear logic.
\newblock In {\em Category Theory and Computer Science}, pages 341--356.
  Springer Berlin Heidelberg, 1989.

\bibitem[dP91a]{dePaiva1991multirelations}
Valeria de~Paiva.
\newblock Categorical {M}ultirelations, {L}inear {L}ogic and {P}etri {N}ets,
  {T}echnical {R}eport, 1991.

\bibitem[dP91b]{dePaiva1991}
Valeria de~Paiva.
\newblock The {D}ialectica {C}ategories.
\newblock University of Cambridge, Computer Lab Technical Report, {PhD} thesis,
  1991.

\bibitem[dP02]{dePaiva2002}
Valeria de~Paiva.
\newblock Lineales: Algebras and categories in the semantics of linear logic.
\newblock {\em Words, Proofs and Diagrams, CSLI Publications, Stanford}, pages
  123--142, 2002.

\bibitem[EDLR16]{emzivat2016probabilistic}
Yrvann Emzivat, Benoit Delahaye, Didier Lime, and Olivier~H Roux.
\newblock Probabilistic time {P}etri nets.
\newblock In {\em International Conference on Application and Theory of Petri
  Nets and Concurrency}, pages 261--280. Springer, 2016.

\bibitem[EFL{\etalchar{+}}20]{Jost2020}
Marzieh Eidi, Amirhossein Farzam, Wilmer Leal, Areejit Samal, and J{\"u}rgen
  Jost.
\newblock Edge-based analysis of networks: curvatures of graphs and
  hypergraphs.
\newblock {\em Theory in Biosciences}, 139(4):337--348, 2020.

\bibitem[EL10]{EHRHARD}
Thomas Ehrhard and Olivier Laurent.
\newblock Interpreting a finitary pi-calculus in differential interaction nets.
\newblock {\em Information and Computation}, 208(6):606--633, 2010.
\newblock Special Issue: 18th International Conference on Concurrency Theory
  (CONCUR 2007).

\bibitem[EW90]{engberg1990petri}
Uffe Engberg and Glynn Winskel.
\newblock Petri nets as models of linear logic.
\newblock In {\em CAAP'90: 15th Colloquium on Trees in Algebra and Programming
  Copenhagen, Denmark, May 15--18, 1990 Proceedings 15}, pages 147--161.
  Springer, 1990.

\bibitem[Gir87]{girard1987}
Jean-Yves Girard.
\newblock Linear logic.
\newblock {\em Theoretical computer science}, 50(1):1--101, 1987.

\bibitem[Ham09]{HAMMACK20091114}
Richard~H. Hammack.
\newblock Proof of a conjecture concerning the direct product of bipartite
  graphs.
\newblock {\em European Journal of Combinatorics}, 30(5):1114--1118, 2009.
\newblock Part Special Issue on Metric Graph Theory.

\bibitem[HLE{\etalchar{+}}09]{synechococcus}
Thomas Hinze, Thorsten Lenser, Gabi Escuela, Ines Heiland, and Stefan Schuster.
\newblock Modelling signalling networks with incomplete information about
  protein activation states: a p system framework of the kaiabc oscillator.
\newblock In {\em International Workshop on Membrane Computing}, pages
  316--334. Springer, 2009.

\bibitem[Hof09]{Hoffmann2009}
Reinhard~W. Hoffmann.
\newblock Introduction.
\newblock In {\em Elements of Synthesis Planning}. Springer Berlin Heidelberg,
  Berlin, Heidelberg, 2009.

\bibitem[H{\"o}h95]{Hohle1995}
U.~H{\"o}hle.
\newblock {\em Commutative, residuated 1---monoids}, pages 53--106.
\newblock Springer Netherlands, Dordrecht, 1995.

\bibitem[IK00]{Imrich2000}
Wilfried Imrich and Sandi Klav{\v z}ar.
\newblock {\em Product graphs, structure and recognition}.
\newblock Wiley, New York, 2000.

\bibitem[Koc20]{kock2020elements}
Joachim Kock.
\newblock Elements of {P}etri nets and processes.
\newblock {\em arXiv preprint arXiv:2005.05108}, 2020.

\bibitem[KSW97]{katis1997nets}
Piergiulio Katis, Nicoletta Sabadini, and Robert~FC Walters.
\newblock Representing place/transition nets in span(graph).
\newblock In {\em International Conference on Algebraic Methodology and
  Software Technology}, pages 322--336. Springer, 1997.

\bibitem[LRSJ21]{forman2018}
Wilmer Leal, Guillermo Restrepo, Peter~F. Stadler, and J{\"{u}}rgen Jost.
\newblock Forman-ricci curvature for hypergraphs.
\newblock {\em Advances in Complex Systems}, 24:2150003, 2021.

\bibitem[Mac71]{maclane}
Saunders MacLane.
\newblock {\em Categories for the Working Mathematician}.
\newblock Springer-Verlag, New York, 1971.
\newblock Graduate Texts in Mathematics, Vol. 5.

\bibitem[Mas19]{master2019}
Jade Master.
\newblock Generalized {P}etri nets.
\newblock {\em arXiv preprint arXiv:1904.09091}, 2019.

\bibitem[Mas20]{master2020petri}
Jade Master.
\newblock Petri nets based on {L}awvere theories.
\newblock {\em Mathematical Structures in Computer Science}, 30(7):833--864,
  2020.

\bibitem[MM90]{meseguer1990}
Jos{\'e} Meseguer and Ugo Montanari.
\newblock {P}etri nets are monoids.
\newblock {\em Information and computation}, 88(2):105--155, 1990.

\bibitem[MOM05]{marti2005petri}
Narciso Marti-Oliet and Jos{\'e} Meseguer.
\newblock From petri nets to linear logic.
\newblock In {\em Category Theory and Computer Science: Manchester, UK,
  September 5--8, 1989 Proceedings}, pages 313--340. Springer, 2005.

\bibitem[NM20]{PhysRevResearch.2.043135}
Zachary~G. Nicolaou and Adilson~E. Motter.
\newblock Missing links as a source of seemingly variable constants in complex
  reaction networks.
\newblock {\em Physical Review Research}, 2:043135, Oct 2020.

\bibitem[RE96]{rozenberg1996elementary}
Grzegorz Rozenberg and Joost Engelfriet.
\newblock Elementary net systems.
\newblock In {\em Advanced Course on Petri Nets}, pages 12--121. Springer,
  1996.

\bibitem[RSS14]{pawel2014reachability}
Julian Rathke, Pawe{\l} Soboci{\'n}ski, and Owen Stephens.
\newblock Compositional reachability in {P}etri nets.
\newblock In {\em International Workshop on Reachability Problems}, pages
  230--243. Springer, 2014.

\bibitem[Sch98]{Schummer1998}
Joachim Schummer.
\newblock The chemical core of chemistry {I}: A conceptual approach.
\newblock {\em Hyle}, 4(2):129--162, 1998.

\bibitem[SNMM11]{Saito2011}
Ayumu Saito, Masao Nagasaki, Hiroshi Matsuno, and Satoru Miyano.
\newblock {\em Hybrid Functional {P}etri Net with Extension for Dynamic Pathway
  Modeling}, pages 101--120.
\newblock Springer London, London, 2011.

\bibitem[SS18]{Stadler2018}
B{\"a}rbel M.~R. Stadler and Peter Stadler.
\newblock Reachability, connectivity, and proximity in chemical spaces.
\newblock {\em MATCH Communications in Mathematical and in Computer Chemistry},
  80(3):639--659, 2018.

\bibitem[ST14]{Seal_Tholen_2014}
Gavin~J. Seal and Walter Tholen.
\newblock {\em Monoidal structures}, page 18–144.
\newblock Encyclopedia of Mathematics and its Applications. Cambridge
  University Press, 2014.

\bibitem[Win87]{Winskel1987}
Glynn Winskel.
\newblock Petri nets, algebras, morphisms, and compositionality.
\newblock {\em Information and Computation}, 72(3):197--238, 1987.

\bibitem[Win88]{winskel1988}
Glynn Winskel.
\newblock A category of labelled {P}etri nets and compositional proof system.
\newblock In {\em [1988] Proceedings. Third Annual Symposium on Logic in
  Computer Science}, pages 142--154. IEEE, 1988.

\end{thebibliography}
\end{document}